\documentclass[12pt, amsfonts, epsfig]{amsart}
\usepackage{latexsym, amssymb, amsmath, epsfig, amscd, amsfonts}
 \newlength{\baseunit}               
 \newcount{\numlines}                
\newcommand{\Z}{\mathbb{Z}}

\newcommand{\R}{\mathbb{R}}
\newcommand{\C}{\mathbb{C}}

        \newfont{\hollow}{msbm10 scaled\magstep1}
        
        \newfont{\Bfmit}{eufm10 scaled\magstep1}
        \newcommand{\bfmit}[1]{\hbox{\Bfmit {#1}}}

\newcommand{\cC}{{\mathcal C}}

\newcommand{\PGL}{\operatorname{PGL}}
\newcommand{\GL}{\operatorname{GL}}

\newcommand{\id}{\operatorname{id}}

\newcommand{\Gr}{\operatorname{Gr}}

\newcommand{\ev}{\operatorname{ev}}

\newcommand{\Proj}{\operatorname{Proj}}

\newtheorem{thm}{Theorem}[section]

\newtheorem{lem}[thm]{Lemma}

\newtheorem{prop}[thm]{Proposition}

\newtheorem{assum}[thm]{Assumption}

\theoremstyle{definition}
\newtheorem{defn}[thm]{Definition}

\newtheorem{exmp}[thm]{Example}

\newtheorem{rem}[thm]{Remark}           

\theoremstyle{remark}

\newcommand{\lremind}[1]{{}}
\newcommand{\bremind}[1]{{}}

\newcommand{\cut}[1]{}

\begin{document}
\pagestyle{plain} \title{{ \large{Topological Aspects of Chow
Quotients} } }
\author{Yi Hu}

\address{Department of Mathematics, University of Arizona, Tucson,
AZ 86721, USA}

\email{yhu@math.arizona.edu}

\address{Center for Combinatorics, LPMC, Nankai University, Tianjin 300071, China}

\maketitle


{\parskip=12pt 

\bigskip
\section{Introduction}

Moduli spaces of geometric objects (e.g., vector bundles,
algebraic curves, etc.) have played central roles in many theories
in algebraic geometry and in its neighboring fields. The
constructions of moduli spaces are frequently done by expressing
them as quotients of schemes by reductive algebraic groups.

However, taking quotients in algebraic geometry is much subtler
than it may appear. Mumford, based upon Hilbert's invariant
theory, developed a systematic method, the Geometric Invariant
Theory (\cite{GIT}),  to deal with projective quotients. There are
several other quotient theories, among them are \cite{Ko97} and
\cite{KM} which construct quotients as algebraic spaces.

It has become well known now that, for a reductive algebraic group
action on a smooth projective variety, Mumford's quotients depend,
in a flip-flop fashion, on choices of linearized line bundles
(\cite{DH} and \cite{Th}). Nevertheless, it is a drawback that
none of Mumford's quotients  is  in general canonical. Besides
this, the closures of the orbits parameterized by a GIT quotient
almost always belong to different cohomology classes, and this,
among other reasons,
 oftentimes makes GIT quotients rather misbehaved compactifications. This is
unsatisfactory from the viewpoint of  moduli problem, where a
moduli space  always parameterizes geometric objects of same
topological type, and awkward to use for purpose of some geometric
computations. To overcome these drawbacks, we are led to consider
a canonical quotient, the Chow quotient. There is another
canonical quotient, the Hilbert quotient. Despite the fact that
the Hilbert quotient, derived from a Hilbert scheme, enjoys more
functorial properties, the Chow quotient, as it parameterizes
cycles, is more geometrically friendly and approachable.

The Chow quotient, $X/\!/^{ch}G$,  of a projective variety $X$ by a  reductive group $G$
is introduced by Kapranov-Sturmfels-Zelevinsky \cite{KSZ} for
toric varieties and by Kapranov \cite{Ka} in general. The
definition of a general Chow quotient is very obscure --- it is
defined as the closure of some Zariski open subset in an ambient
compact variety (Chow variety, to be precise, another obscure
space).

The main purpose of this paper is to give some topological
interpretations and characterization of Chow quotient which have
the advantage to be  more intuitive and geometric. This is to be
done over the field of complex numbers and  in the languages that
are familiar to topologists and differential geometers. Here, we
are content to focus on torus actions even though some of our
results remain true for more general reductive groups. (The case
of general group actions will be treated elsewhere.)

To give the reader some good ideas about Chow quotient as well as
its various topological characterization and interpretations that
we will formally introduce in the main body of the paper, let us
{\it informally} consider a simple, yet quite informative example.

Let $G = \C^*$ act on ${\mathbb P}^2$ by $$ \lambda \cdot [x:y:z]
= [\lambda x: \lambda^{-1} y : z].$$ Consider a map $\Phi:
{\mathbb P}^2 \rightarrow \R$ given by $$\Phi([x:y:z]) =
\frac{|x|^2 - |y|^2}{|x|^2 + |y|^2 + |z|^2}.$$ This is the moment
map for the induced symplectic $S^1$-action with respect to the
Fubibi-Study metric. Its image is the interval $[-1,1]$.

\bigskip

\begin{picture}(20, 9)
\put(4,.2){ \psfig{figure=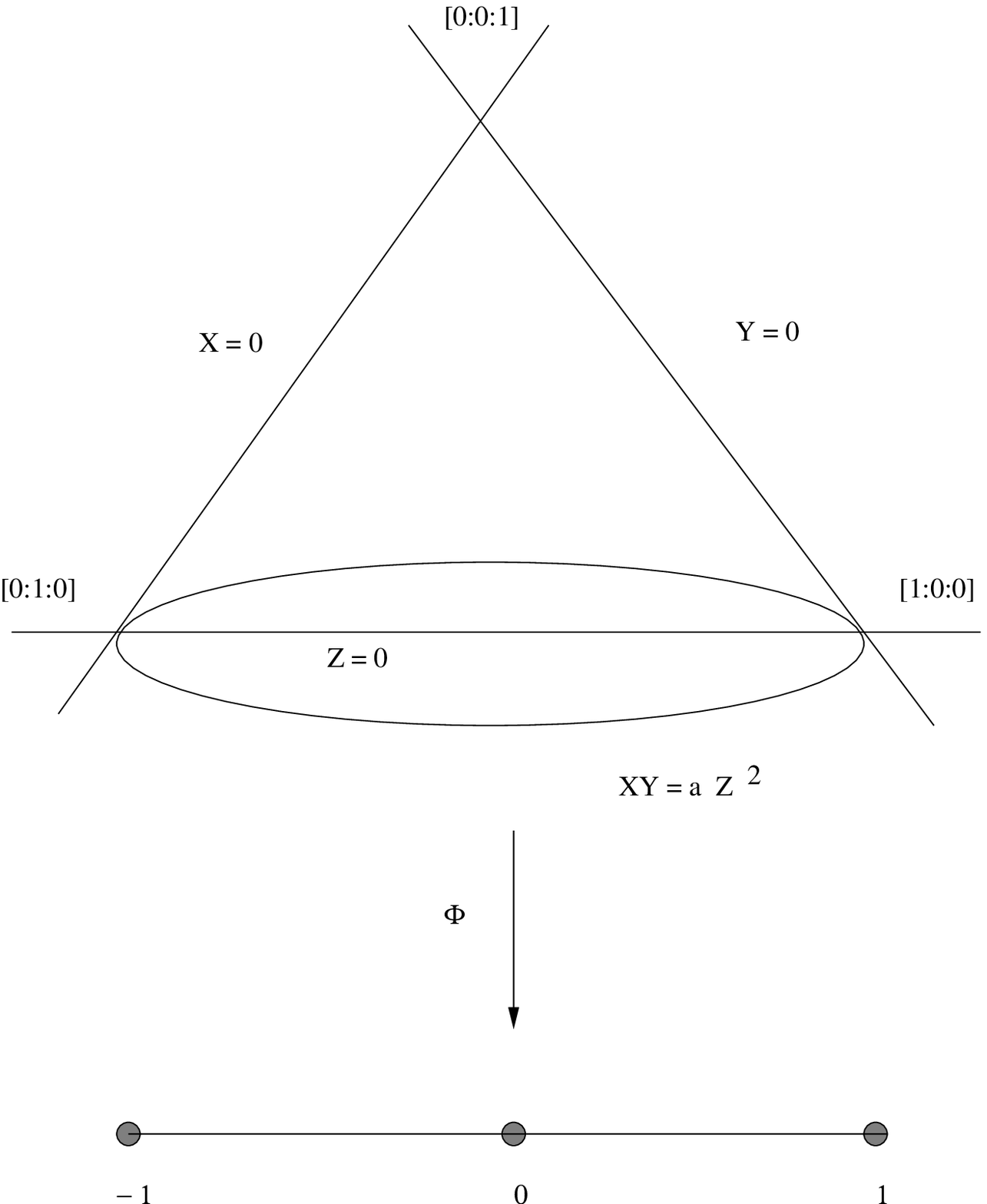, height = 8.5cm,width = 8cm}
}
\end{picture}

\bigskip
\centerline{Figure 1. \ Conics}





The $\C^*$ orbits are classified as follows.  (See Figure 1 for an
illustration.)
\begin{itemize}
\item Generic $\C^*$-orbits are conics $XY = a Z^2$ minus two
points [1:0:0] and [0:1:0] for $a \ne 0, \infty$. (In this
introduction, we will always assume that $a \ne 0, \infty$.) We
denote these orbits by $(XY = a Z^2)$. The moment map image of the
orbit $(XY = a Z^2)$ is $(-1,1)$.
 \item Other 1-dimensional orbits
are the three coordinate lines $X=0$, $Y=0$, and $Z=0$ minus the
coordinate points on them. We denote these orbits by $(X=0)$,
$(Y=0)$, and $(Z=0)$. The moment map images of the orbits $(X=0)$,
$(Y=0)$, and $(Z=0)$ are $(-1,0), (0,1)$, and $(-1,1)$,
respectively. \item Finally, the fixed points are the three
coordinate points, [1:0:0], [0:1:0], and [0:0:1], and their moment
map images are $1, -1,$ and $0$, respectively.
\end{itemize}
 The moment
map has three critical values $-1, 0,$ and $1$ which divide the
interval into two top chambers $[-1,0]$ and $[0,1]$, and three
0-dimensional chambers $\{-1\}, \{0\}, \{1\}$. Each chamber $C$
defines a GIT stability:  a point $[x:y:z]$ is semi-stable with
respect to the chamber $C$ if $C \subset \Phi(\overline{\C^* \cdot
[x:y:z]})$, and it is stable if the (relative) interior $C^\circ
\subset \Phi(\C^* \cdot [x:y:z])$ and $\dim \C^* \cdot [x:y:z] =
1$. (For a reference for this, see for example, \cite{Hu91}.)
Thus, the orbit $(X=0)$  is stable with respect to $[-1,0]$,
unstable with respect to $[0,1]$; while the orbit $(Y=0)$ is
stable with respect to $[1,0]$, unstable with respect to $[-1,0]$.
But, $(X=0)$, $(Y=0)$, and $[0:0:1]$ are all semi-stable with
respect to the chamber $\{0\}$.

\bigskip

\begin{picture}(12, 11)
\put(3,1){ \psfig{figure=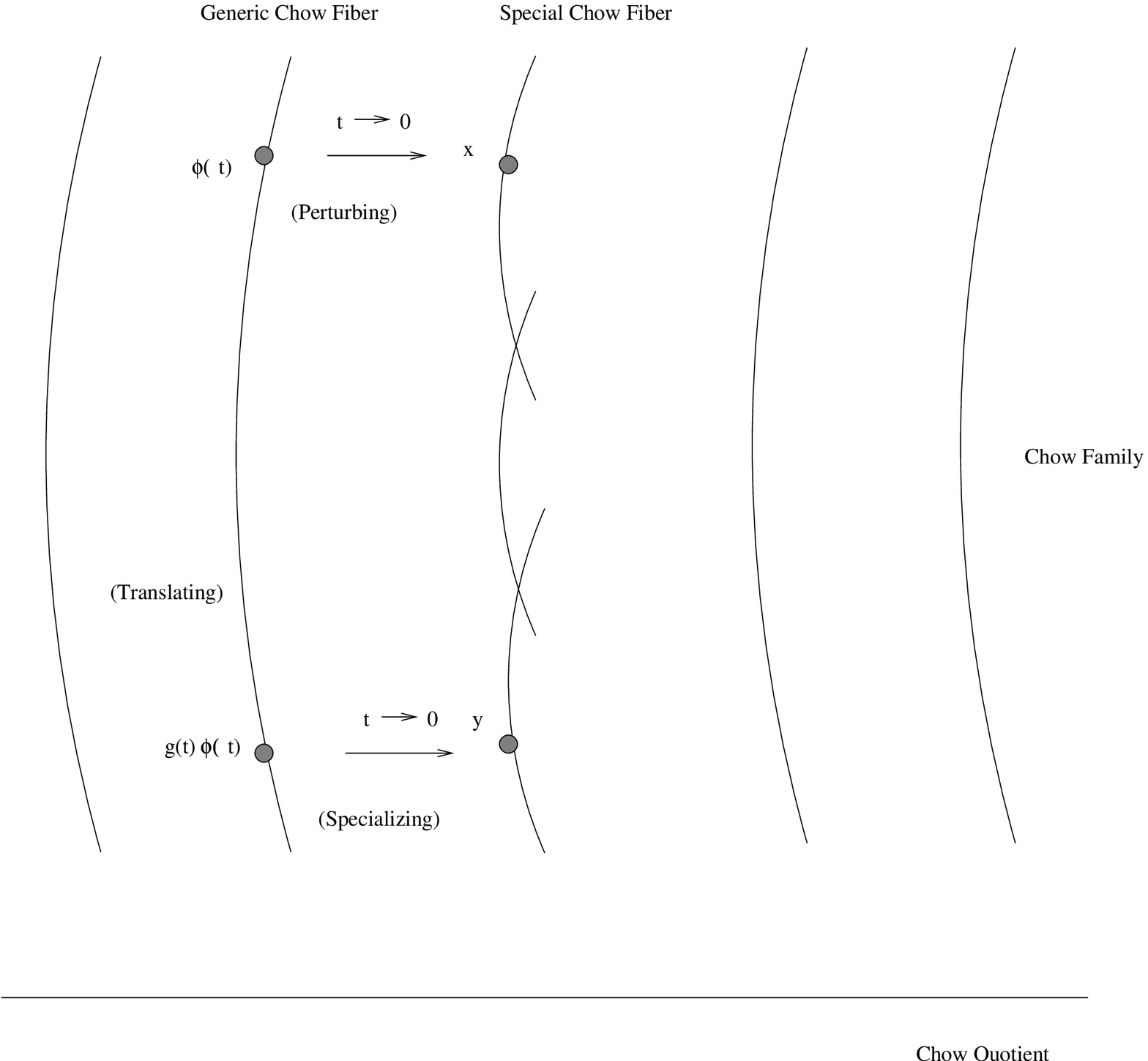,height = 9.5cm,width = 10.5cm} }
\end{picture}

\centerline{Figure 2.  Perturbing-Translating-Specializing}

A general feature of GIT quotient is that it parameterizes orbits
that are {\it closed} in the semi-stable locus. Hence, the GIT
quotient $X_{[-1,0]}$ defined by the chamber $[-1,0]$
parameterizes $(XY = a Z^2)$, $(Z=0)$, and $(X=0)$. The GIT
quotient $X_{[0,1]}$ defined by the chamber $[0,1]$ parameterizes
$(XY = a Z^2)$, $(Z=0)$, and $(Y=0)$. And, the GIT quotient
$X_{\{0\}}$ defined by the chamber $\{0\}$ parameterizes $(XY = a
Z^2)$, $(Z=0)$, and $[0:0:1]$. Now observe that the lines $X=0$
and $Y=0$, which are of degree 1,  have different homology classes
than the conic orbits $XY = a Z^2$, which are of degree 2. This is
undesirable. Even worst, the two orbits $(X=0)$ and $(Y=0)$
 are all identified with the closed orbit $[0:0:1]$ in the quotient $X_{\{0\}}$,
  and the unique
closed orbit $[0:0:1]$ even has different dimension than those of
the conic orbits.

\bigskip

However, Chow quotient takes completely different approach. It
first considers  the closures of the generic $\C^*$-orbits,  $XY =
a Z^2$ ($a \ne 0, \infty$), and then looks at all their possible
degenerations. When $a=0$, we get the degenerated conic $XY=0$,
two crossing lines; and when $a = \infty$, we obtain $Z^2 =0$, a
double line. They all have the same homology classes (degree 2).
And the Chow quotient is {\it the space of all $\C^*$-invariant
conics}. (See Figure 1.) Each point of the Chow quotient
represents an invariant algebraic cycle. In this case, these
cycles are $[XY = a Z^2]$ ($a \ne 0, \infty$), $[X=0] + [Y=0]$,
and $2 [Z=0]$.

From this  example, we see that Chow quotient parameterizes cycles
of generic orbit closures and their toric degenerations which are
{\it certain sums of orbit closures of dimension $\dim G$}. We
will call these cycles {\it Chow fibers}. So, when do two
arbitrary points belong to the same Chow fiber? Consider the
example again. We have that $[0:y:z]$ and $[x:0:z]$ ($xyz \ne 0$)
belong to the same Chow fiber $XY=0$. To get $[x:0:z]$ from
$[0:y:z]$, we first {\it perturb} $[0:y:z]$ to a general position
$\varphi (t) = [tx:y:z]$ ($t \ne 0$),  then {\it translate} it by
$g(t) = t^{-1} \in \C^* $ to $g (t) \cdot \varphi (t) = [x:ty:z]$,
and then $g (t) \cdot \varphi (t)$ {\it specializes} to $[x:0:z]$
when $t$ goes to 0. We will call this process {\it
perturbing-translating-specializing} (PTS). It turns out this
simple relation holds true in general. That is, we prove in
general that two points $x$ and $y$ of $X$, with $\dim G \cdot x =
\dim G \cdot y = \dim G$, belong to the same Chow fiber if and
only if $x$ can be perturbed (to general positions), translated
along $G$-orbits (to positions close to $y$), and then specialized
to the point $y$. See Figure 2 for an illustration. This is our
Theorem \ref{chowfiberAG}.

An upshot of PTS relation is that, comparing to the nondescriptive
definition of special Chow fibers, it is computable, and thus
provides some much needed information on boundary cycles of the
Chow quotient. As an application, we apply Theorem
\ref{chowfiberAG} to the case of point configurations on ${\Bbb
P}^n$ ($n >1$), and propose a {\it geometric interpretation} of
the Chow quotients of $({\Bbb P}^n)^m$ (equivalently, the Chow
quotients of higher Grassmannians). It was relayed to me that
Lafforgue's space constructed in \cite{Lafforgue} is in general
reducible, one component of which is the Chow quotient, and,
Theorem \ref{chowfiberAG} may be used to give a way of telling
Lafforgue's components apart. (We will pursue the above and
related issues in a forthcoming paper.)

Back to our example, let $$G= \C^* = S^1 \times {\mathbb R}_{>0} =
K \times A$$ be the polar decomposition. Then, one checks easily
that the moment map $\Phi$ is $S^1$-invariant and $\Phi( r \cdot
[x:y:z])$ is a strictly increasing function of $r$ for $r \in
{\mathbb R}_{>0}$ and for any given non-fixed point $[x:y:z]$.
This implies, for example,  that in each closed orbit
parameterized by a point of the GIT quotient $X_{[-1,0]}$ there is
a unique $S^1$-orbit whose moment map image is, say, $-
\frac{1}{2}$ (any other point in $(-1,0)$ will work equally well).
This in turn implies that $X_{[-1,0]}$ is homeomorphic to
$\Phi^{-1}(- \frac{1}{2})/S^1$.  Similarly, we have $X_\{0\} \cong
\Phi^{-1}(0)/S^1$ and $X_{[0,1]} \cong
\Phi^{-1}(\frac{1}{2})/S^1$. In particular, $\Phi^{-1}(-
\frac{1}{2})/S^1$ has a unique point parameterizing the orbit
$(X=0)$,  $\Phi^{-1}(\frac{1}{2})/S^1$ instead has a point
parameterizing the orbit $(Y=0)$, while $\Phi^{-1}(0)/S^1$ has a
point parameterizing the smaller orbit $(X=0=Y)$ ($=[0:0:1]$). As
it turns out, these are just examples of the very general
correspondence between GIT quotients and symplectic quotients. Put
it more formally (cf. \cite{Kirwan}, \cite{GIT}), symplectic
quotients endow various symplectic structures, possibly singular,
on GIT quotients. Pushing this circle of ideas further, it is
natural to ask whether a Chow quotient admits its own symplectic
counterparts.

\bigskip

\begin{picture}(15, 9)
\put(4,1){ \psfig{figure=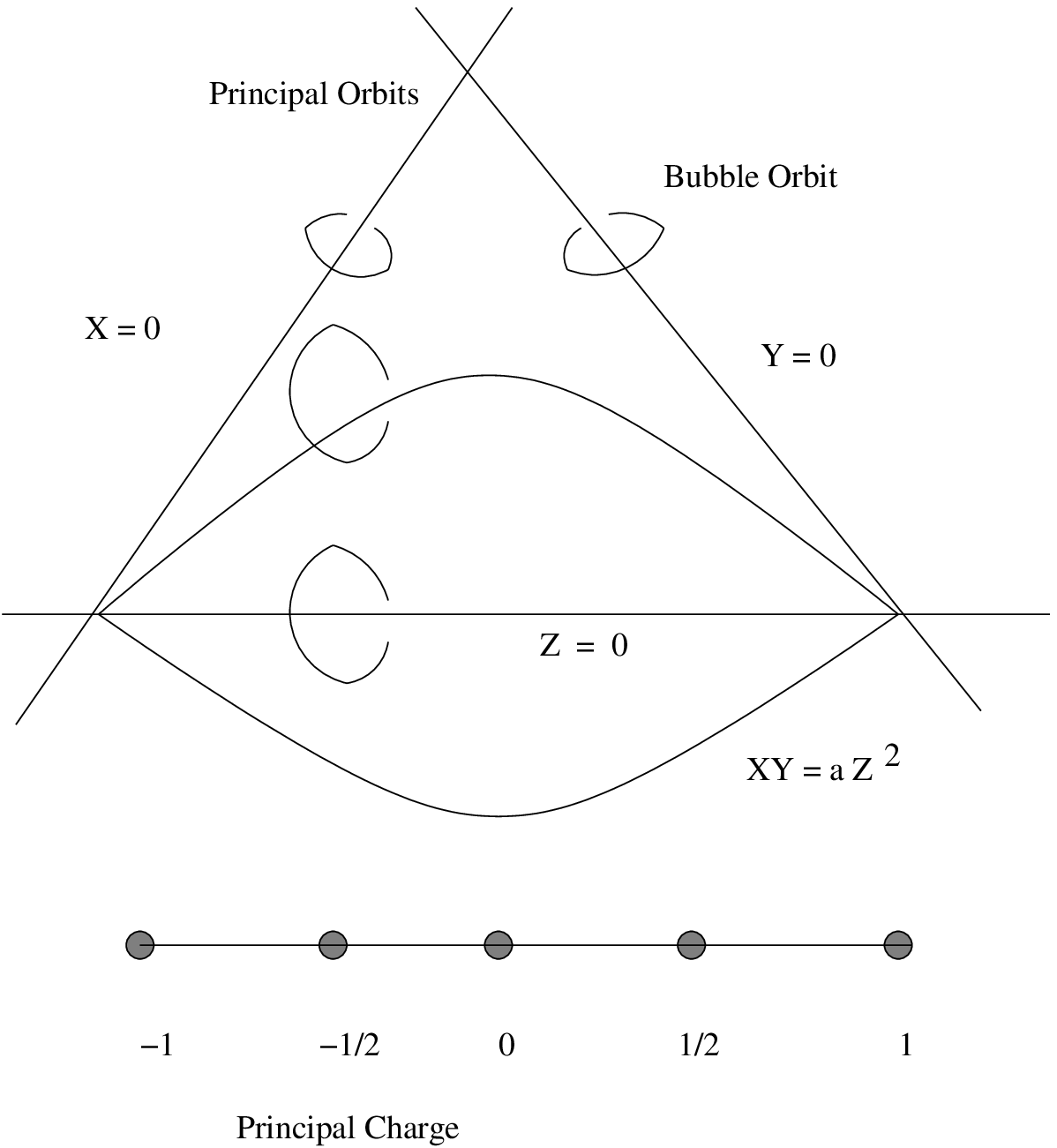, height = 7.5cm,width =
7cm} }
\end{picture}

\centerline{Figure 3.  Stable $S^1$-Orbits}

To this end, we introduce the so-called {\it stable $K$-orbits}
with a fixed set  of {\it momentum charges}, $\gamma$. Consider
again the Chow quotient of our previous example, the space of all
invariant conics in ${\mathbb P}^2$. So, here  $K$= $S^1$.  We
first fix a general point, say $-\frac{1}{2}$, in $[-1, 1]$, and
call it the {\it principal momentum charge}.  Then each generic
conic orbit contains a unique $S^1$-orbit whose moment map image
is $-\frac{1}{2}$. These $S^1$-orbits are the generic stable
$K$-orbits with the principal momentum charge $-\frac{1}{2}$. The
same can be done for the orbit $(Z=0)$. That is, there is a unique
$S^1$-orbit in $(Z=0)$ whose moment map image is principal
momentum charge $-\frac{1}{2}$. The point here is that
$\Phi(XY=aZ^2) = \Phi (Z=0) = [-1,1]$ ($a \ne 0, \infty$), hence
$\Phi(XY=aZ^2)$ and  $\Phi (Z=0)$ do not subdivide $[-1,1]$. For
the cycle $[X=0] + [Y=0]$, the moment map image $\Phi (X=0)$ and
$\Phi (Y=0)$ form a subdivision of $[-1, 1] = [-1,0] \cup [0, 1]$.
In this case, in $[-1, 0]$ we (have to) stick with the principal
momentum charge $-\frac{1}{2}$, but in $[0,1]$ we may choose any
interior point, say, $\frac{1}{2}$. Then we have two unique
$S^1$-orbits $S^1 \cdot [0:1:1]$ and $S^1 \cdot [1:0:1]$ in the
cycle $[X=0] + [Y=0]$ over the prescribed momentum charges
$\{-\frac{1}{2}, \frac{1}{2}\}$. The union of these two orbits
$S^1 \cdot [0:1:1] \cup S^1 \cdot [1:0:1]$ is what we call a
special stable $S^1$-orbit with the prescribed set of  momentum
charges. (See Figure 3\footnote{A few words on Figure 3. Note that
each $G$-orbit closure is a piecewise fibration over the moment
map image with the fibers the compact group orbits. In Figure 3,
we depict the circle fibers over the principal charge
$-\frac{1}{2}$ and the charge $\frac{1}{2}$. A circle winding
around a curve or a line means that the circle is contained in the
cycle represented by the curve or the line.}.) Then, we see that
the Chow quotient also parameterizes stable $S^1$-orbits with a
fixed set of momentum charges $\gamma$ (in this example, $\gamma$
consists of $-\frac{1}{2}$ and $\frac{1}{2}$). We put all stable
$S^1$-orbits together, denote it by ${\bfmit M}_{\gamma}$, and
call it  the moduli space  of stable $K$-orbits with prescribed
momentum charges $\gamma$. We prove in general (see Assumption
\ref{uniquesupport}) that the space ${\bfmit M}_\gamma$ is always
homeomorphic to the Chow quotient $X/\!/^{ch}G$,  regardless of
the choice of $\gamma$ (Theorem \ref{moduliofstabelKorbits}).

Now,  we insert here a favorable property enjoyed by the Chow
quotient: it dominates all GIT quotients. This was proved by
Kapranov, and the dominating morphisms were also discovered in
\cite{Hu91}. In fact, the Chow quotient , under Assumption
\ref{uniquesupport}, is the {\it least common refinement} of all
GIT quotients, in a {\it strict} sense. That is, $X/\!/^{ch}G$ is
homeomorphic to the limit quotient $X/\!/^{lim}G$,  the
distinguished irreducible component of the inverse limit of all
GIT quotient (Theorem \ref{chow=limit}). When $X$ is a toric
variety, this is proved in \cite{KSZ}.

Then, the morphism from the Chow quotient to a GIT quotient,
 in terms of stable $K$-orbits, corresponds to a map from
 ${\bfmit M}_\gamma$ to the symplectic quotient $\Phi^{-1}({\bf r})/K$
 where  ${\bf r}$ is the
principal momentum charge in $\gamma$, and this last map is quite
transparent. Every stable orbit contains a unique principal orbit,
the one with the principal momentum charge, all the rest will be
called ``bubble'' orbits of the principal one. Then, the map from
${\bfmit M}_\gamma$ to $\Phi^{-1}({\bf r})/K$ simply forgets all
the bubble orbits and send a stable orbit to its principal part.
See Figure 3.

The whole picture of the stable $K$-orbits resembles that of
stable polygons in \cite{Hu99a} where we give a symplectic
construction of $\overline{M}_{0,n}$, the moduli space of stable
$n$-pointed rational curves. In fact, stable polygon is the source
of inspiration for the introduction of stable $K$-orbit, and,
\cite{Hu99a} may be viewed as an (interesting and long) example of
the general theory described in this paper.

The space ${\bfmit M}_\gamma$ reflects the Hamiltonian aspects of
the $G$-action on $X$. Transversal to Hamiltonian flows are the
gradient flows of various 1-dimensional projections of the moment
map. Indeed, let $G = K \cdot A$ be the polar decomposition of
$G$. Then Hamiltonian flows are tangential to $K$-orbits, while
gradient flows are tangential to $A$-orbits. This shift of
viewpoint leads us to find another topological approach to the
Chow quotient.   Here, to make things work coherently, we
introduce  the so-called {\it stable action-manifolds}: generic
stable action-manifolds are simply  the closures of $A$-orbits
through generic points;
 special ones are certain configurations of
closures of $A$-orbits that are resulted as the limits of generic
ones.  In the example, generic stable action-manifolds are
$\overline{\R_{>} \cdot [a:1:1]}$ ($a \in \R, \ne 0, \infty$). The
special action-manifolds are $\overline{\R_{>} \cdot [1:1:0]}$ and
$\overline{\R_{>} \cdot [0:1:1]} \cup \overline{\R_{>} \cdot
[1:0:1]}$. Up to $K$-action, these are all the piece-wise gradient
flow lines of the map $\Phi$ connecting $[1:0:0]$ and $[0:1:0]$.
See Figure 4. In general, stable action-manifolds are piece-wise
smooth manifolds with corners.

\bigskip

\begin{picture}(6, 7)
\put(4,1){ \psfig{figure=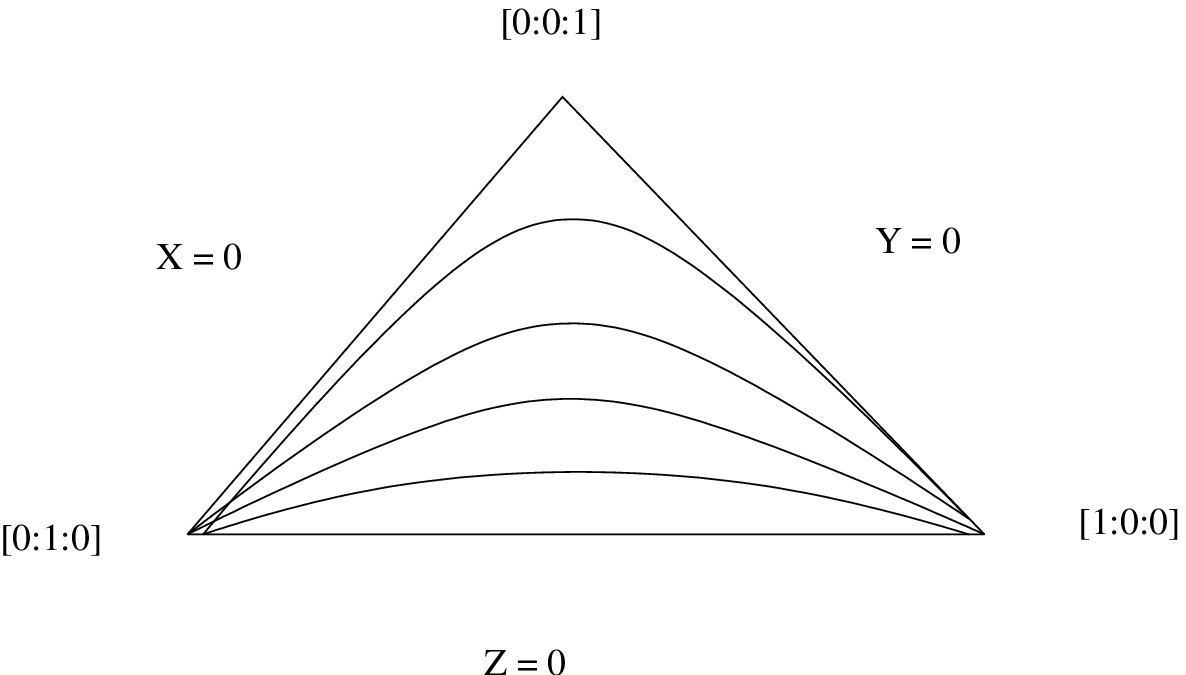, height = 6cm,width = 7cm}
}
\end{picture}

\centerline{Figure 4.  Stable Action-Manifolds}

 To rigorously construct the special
 configurations,  we apply a real version of our theorem on PTS,
 which is formulated (for the ad hoc purpose), and proved
in \S \ref{section:toppush}. Two stable action-manifolds are said
to be equivalent if one can be transferred to the other by the
action of an element of the compact group $K$. Then we prove,
again under Assumption \ref{uniquesupport}, that the moduli space
${\bfmit M}$ of equivalence classes of stable action-manifolds
exists as a separated complex variety, and is homeomorphic to the
Chow quotient. Stable action-manifolds are orthogonal to that of
stable $K$-orbits. Hence in this sense, the space ${\bfmit
M}_\gamma$ parameterizes Hamiltonian slices  of Chow fibers, while
the moduli of stable action manifolds ${\bfmit M}$ parameterizes
gradient slices (with respect to the moment map) of Chow fibers
which are orthogonal to Hamiltonian slices.

In the previous example, we worked out all sort of the quotients
considered in this paper (GIT and Chow), but they are all
isomorphic to ${\mathbb P}^1$, the unique compactification of
$\C^*$,  because we insist an example that are very simple to
describe. Here it should be fair to at least point out to the
reader a workable example where a nontrivial wall crossing
phenomenon and different quotients do occur. Hence we take the
liberty to include the following example with  details left to the
reader. So, consider the action of $\C^*$ on ${\mathbb P}^3$ by
$$ \lambda \cdot [x:y:z:w] = [\lambda x: \lambda y : \lambda^{-1}z: w].$$
The moment map is
$$\Phi([x:y:z:w]) = \frac{|x|^2 + |y|^2 - |z|^2}{|x|^2 + |y|^2 + |z|^2 + |w|^2}.$$

The image $\Phi(X)$ is $[-1,1]$ with three critical values $-1, 0,
1$.  So, we consider the level sets $\Phi^{-1}(-\frac{1}{2}),
\Phi^{-1}(0), \Phi^{-1}(\frac{1}{2})$. In Figure 5, we illustrate
the {\it real parts} $\Phi^{-1}(-\frac{1}{2})_\R, \Phi^{-1}(0)_\R,
\Phi^{-1}(\frac{1}{2})_\R$ of the level sets restricted to $\C^3
\subset {\mathbb P}^3$ (the $\C^3$ is defined by setting $w= 1$).
It turns out this real picture preserves all the topological
information we need.

To understand this picture, note that the real points of $S^1$ are
$\Z_2 = \{-1, 1\}$. So, the real parts of
$$\Phi^{-1}(-\frac{1}{2})/S^1, \Phi^{-1}(0)/S^1,
\Phi^{-1}(\frac{1}{2})/S^1$$  are
$$\Phi^{-1}(-\frac{1}{2})_\R/\Z_2, \Phi^{-1}(0)_\R/\Z_2,
\Phi^{-1}(\frac{1}{2})_\R/\Z_2.$$

\begin{picture}(15, 9)
\put(4,1){ \psfig{figure=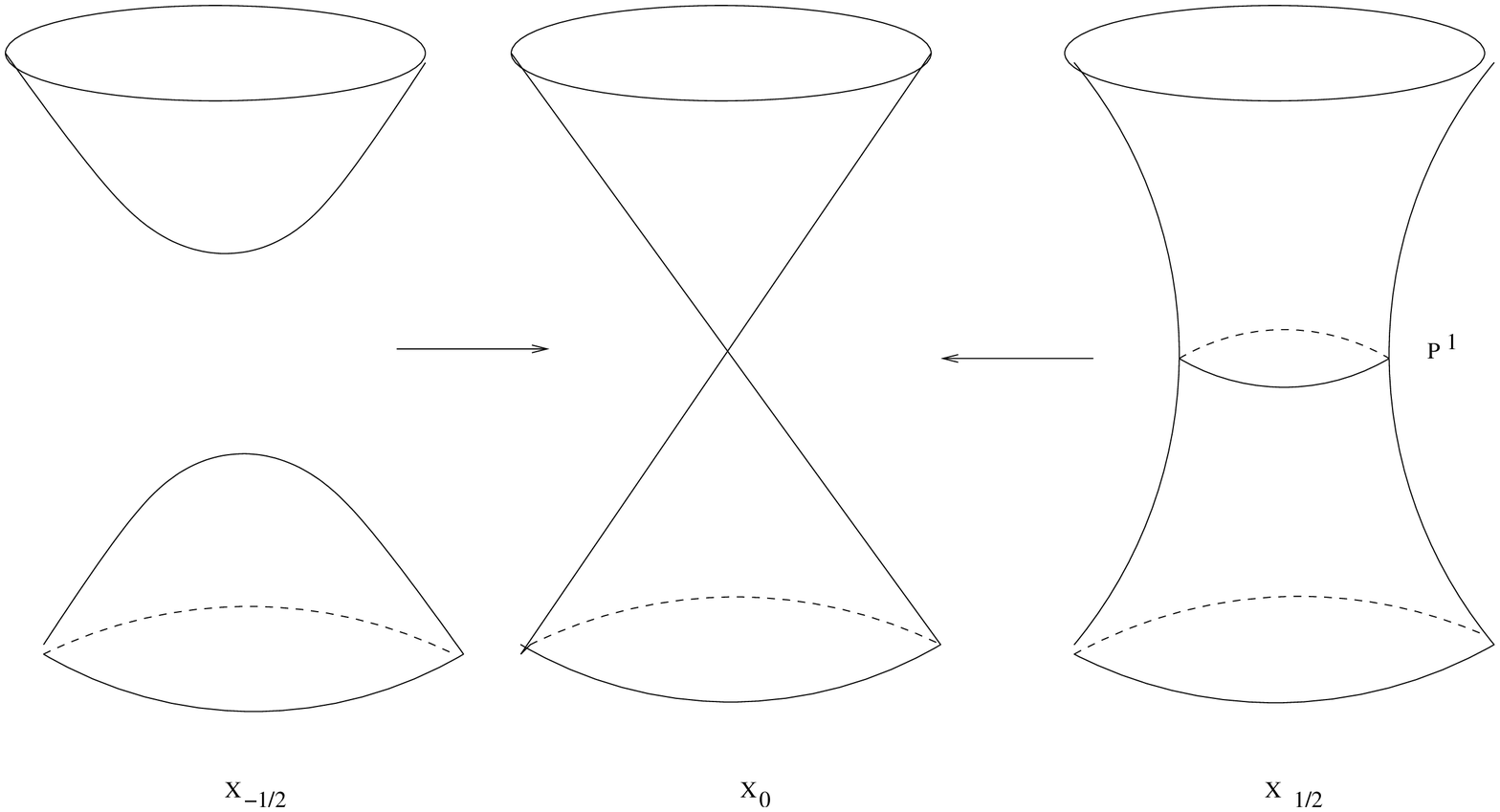,height = 7.5cm,width = 8cm} }
\end{picture}

\centerline{Figure 5.  Wall-Crossing Maps}

\bigskip

Note that $\Z_2$ acts on each level set by identifying the lower
part with the upper part. Hence, the quotient can be naturally
identified with the upper part. Now observe that the left map
happens to be an isomorphism, but the right map is a (real) blowup
along the origin so that the special fiber is $\R {\mathbb P}^1$.
Now complexifying this picture and compactifying the results, we
obtain three quotients: $X_{[-1,0]} \cong {\mathbb P}^2$,
$X_{\{0\}} \cong {\mathbb P}^2$, and $X_{[0,1]}$ isomorphic to the
blowup of ${\mathbb P}^2$ along a point. The reader may try to
classify the generic $\C^*$-orbits and study their degenerations.
He can verify that the Chow quotient is also isomorphic to the
blowup of ${\mathbb P}^2$ along a point.

I learned from Igor Dolgachev, during the Arizona conference on
Geometry and Topology of Quotients (Dec. 5-8, 2002),
 that Yuri Neretin first considered the PTS
relation\footnote{Neretin did not use this terminology.  Keel is
partially responsible for our choice of the term.} (\cite{N1},
\cite{N2}), and conjectured that it may relate to some quotient
construction. Soon after the conference, I realized that the
relation actually characterizes Chow fibers (Theorem
\ref{chowfiberAG}), and its real topological modification can be
used to geometrically compactify the moduli space of generic
action-manifolds (Theorem \ref{MofSA}), improving our original
approach. It seems from this paper that Neretin's quotient and
``hinges'' for symmetric spaces (\cite{N1}, \cite{N2}) are analogs
of Chow quotient and Chow fibers in the topological
situation\footnote{ This was observed through a communication with
Dolgachev after posting this paper on ArXiv.}. Some result here
may be expected by Dolgachev and Neretin.

\bigskip\noindent
{\sl Acknowledgments}. I am very grateful to the anonymous referee
who carefully read the earlier versions of this paper and offered
numerous very critical comments and constructive suggestions. As a
result, the exposition of this paper is substantially improved and
hopefully more readable. Financial support and hospitality from
Harvard University and Professor S.-T. Yau, from NCTS Taiwan and
Professor C.L Wang, and from Hong Kong UST and Professors W.-P. Li
and Y. Ruan are gratefully acknowledged. This research was
partially supported by NSF and NSA.

\bigskip



\bigskip
\section{GIT quotient and Symplectic Reduction}

\subsection{GIT quotient}
\label{GIT}
Let the torus $G = {\Bbb C}^n$ act on a smooth projective variety $X$ over the field of complex numbers.
Throughout the paper, we will assume that the action is generically free, that is,
the isotropy subgroups of generic points are the identity subgroup.

Let $K = (S^1)^n$ be the compact form of $G$. Fix an ample
linearized line bundle $L$ on $X$. Now pick a $K$-invariant
Hermitian metric on $L$; equivalently, a $K$-invariant symplectic
form in $[c_1(L)]$. Then, there is an uniquely associated moment
map
$$\Phi_L: X \rightarrow {\bfmit k}^*$$
where ${\bfmit k}^*$ is the linear dual of the Lie algebra ${\bfmit k}$ of $K$.

The moment map image $P_L = \Phi_L(X)$ is a compact polytope
(\cite{Atiyah}, \cite{GS}). Atiyah also shows that $\Phi_L
(\overline{G \cdot x})$ is a subpolytope of $P_L$ for any $x \in
X$, and $\Phi_L (\overline{G \cdot x})= P_L$ for generic $x$.
 For simplicity, we often simply write $P$
instead of $P_L$. This polytope admits a natural chamber
decomposition, $\cC_L$,  by the common refinement of all
the subpolytopes of the form
$$\Phi_L(\overline{G \cdot x}), x \in X.$$
For any point ${\bf r} \in P$, we will use $[{\bf r}]$ to denote the minimal chamber that
${\bf r}$ belongs to.

Given any rational point ${\bf r} \in P$, an integral multiple $m
{\bf r}$ can be identified with a character $\chi$ of $G$. Let
$L^m(\chi)$ be the linearized line bundle $L^m$ twisted by
$-\chi$. In the symplectic terms, this means we replace the moment
map $\Phi_L$ by $\Phi'=m\Phi_L - m{\bf r}$ (this is equivalent to
$\Phi_L - {\bf r}$ as far as stability is concerned). Then we will
obtain a Zariski open subset $X^{ss}(L^m(\chi))$, the semi-stable
locus with respect to the linearization $L^m(\chi)$. By Theorem
8.3 of \cite{GIT}, a point $x \in X$ is semi-stable with respect
to  $L^m(\chi)$ if $0 \in \Phi'(\overline{G \cdot x})$ and it is
stable if $0 \in \Phi'(G \cdot x)$ and $\dim G \cdot x = \dim G$.
It follows that $X^{ss}(L^m(\chi))$ does not depend on the choice
of the multiple $m$,  hence we will denote it by
$$X^{ss}({\bf r}),$$
and call it the set of semi-stable points with respect to ${\bf
r}$. Observe that
$$ x \in X^{ss}({\bf r}) \; \Leftrightarrow \;  {\bf r} \in
\Phi_L(\overline{G \cdot x}).$$
 The GIT quotient $X^{ss}({\bf r})/\!/G$ exists as a separated
projective variety. Topologically, $X^{ss}({\bf r})/\!/G$ is
obtained from $X^{ss}({\bf r})$ by identifying points with the
following equivalence relation:
$$x \in X^{ss}({\bf r}) \sim y \in X^{ss}({\bf r}) \; \hbox{iff}
\; \overline{G \cdot x} \cap \overline{G \cdot y} \cap X^{ss}({\bf
r}) \ne \emptyset.$$ It follows that when all semi-stable points
are actually stable, two different orbits will never be
identified. And, in such a case, $X^{ss}({\bf r})/\!/G$ is
(topologically) the ordinary orbit space.

Moreover, one can easily deduce  that
$$X^{ss}({\bf r})= X^{ss}({\bf r}')$$
if and only if ${\bf r} $ and ${\bf r}'$ are in the relative
interior of the same chamber $C$ for some $C \in \cC_L$. Hence we
may write $X^{ss}(C)$ for $X^{ss}({\bf r})$ for all interior
rational point ${\bf r} \in C$. Observe that in this case
$$ x \in X^{ss}(C) \; \Leftrightarrow \;  C \subset
\Phi_L(\overline{G \cdot x}).$$ In particular, if $C$ is of top
dimension, all semi-stable points are stable.
 To make our notation concise, we
will use $X_{[{\bf r}]}$ or $X_C$ to denote the GIT quotient
$$X^{ss}({\bf r})/\!/G = X^{ss}(C)/\!/G.$$

For any chamber $C$, if $D \subset C$ is a face of $C$, then one
can check from the definition that we have the inclusion
$$X^{ss}(C) \subset X^{ss}(D),$$
and this inclusion induces a canonical birational projective morphism
$$f_{CD}: X_C \rightarrow X_D.$$
They all together form an inverse system
$$\{ X_C,  f_{CD} | D \subset C \in \cC_L\}.$$
For a reference for the above, one may consult \cite{Hu91}.

\subsection{Symplectic reduction}
For any point ${\bf r} \in P$, rational or not, we have an inclusion
$$G \cdot \Phi^{-1}_L({\bf r}) \subset X^{ss}(C),$$
where $C$ is the unique minimal chamber containing ${\bf r} $.
In fact, $G \cdot \Phi^{-1}_L({\bf r})$ is exactly
the set of closed orbits of $X^{ss}(C)$. Hence the above inclusion induces
a natural homeomorphism
$$\Phi^{-1}({\bf r})/K \rightarrow X_C,$$
from the {\it symplectic reduction} $\Phi^{-1}({\bf r})/K$ to the
GIT quotient $X_C$, thanks to a theorem of Kirwan (cf.
\cite{GIT}). This basically says that the GIT quotient $X_C$
carries a family of symplectic structures, possibly singular,
parameterized by the interior points of the chamber $C$.

Again, to be concise,
 we will use $X_{\bf r}$ to denote the symplectic reduction $\Phi^{-1}({\bf r})/K$.
Note that when ${\bf r}$ is rational, we have used $X_{[{\bf r}]}$ to denote the GIT quotient
defined by ${\bf r}$; the  subscript $[{\bf r}]$ emphasizes the fact that the GIT quotient only depends
on the minimal  chamber that ${\bf r}$ belongs to but not the individual point ${\bf r}$.

\bigskip
\section{Chow Quotient: algebro-geometric approach}

\subsection{Chow Variety}
\label{ChowVariety} We recall briefly the Chow variety. The
reference is \cite{GKZ}. For a full account, see \cite{Ko}.

 A $k$-dimensional algebraic cycle in ${\mathbb P}^n$ is a formal
finite linear combination $C = \sum_i m_i C_i$ with non-negative
integer coefficients, where $C_i$ are $k$-dimensional irreducible
closed subvarieties in ${\mathbb P}^n$. The degree of $C$ is
$\sum_i m_i \hbox{deg} (C_i)$. Assume $Y \subset {\mathbb P}^n$ is
a $k$-dimensional irreducible subvariety of degree $d$. Consider
the set $Z(Y)$ of all $(n-k-1)$-dimensional linear subspaces of
${\mathbb P}^n$ that intersect $Y$. Then $Z(Y)$ is an irreducible
hypersurface of degree $d$ in the Grassmannian $\Gr(n-k-1,{\mathbb
P}^n)$. Let
$$B = \oplus_m B_m = \oplus_m H^0(\Gr(n-k-1,{\mathbb P}^n), {\mathcal O}(m))$$
be the coordinate ring of $\Gr(n-k-1,{\mathbb P}^n)$. Then $Z(X)$
is defined by the vanishing of some element $R_Y \in B_d$ which is
unique up to homothety (multiplication by a nonzero scalar). This
element will be called the Chow form of $Y$. Now let $C = \sum_i
m_i C_i$ be an algebraic cycle of degree $d$. We define the Chow
form of $C$ as
$$R_C = \prod_i R_{C_i}^{m_i} \in B_d.$$
Let $\hbox{Chow}_d ({\mathbb P}^n)$ be the set of all
$k$-dimensional algebraic cycles of degree $d$ in ${\mathbb P}^n$.
Then the map
$$\hbox{Chow}_d ({\mathbb P}^n) \rightarrow \Proj(B_d), \;\;\; C \to
R_C$$ defines an embedding of $\hbox{Chow}_d ({\mathbb P}^n)$ into
the projective space $\Proj(B_d)$. The variety $\hbox{Chow}_d
({\mathbb P}^n)$ with the projective algebraic structure defined
by the above embedding is called the Chow variety of
$k$-dimensional algebraic cycles of degree $d$ in ${\mathbb P}^n$.

For a general projective variety $X$, an algebraic cycle is
defined in the same way  as a non-negative linear combination of
irreducible subvarieties. Each irreducible subvariety represents a
homology class. The homology class of an algebraic cycle is
defined as the induced sum of homology classes of the irreducible
components. Now fix a homology class $\delta$. Let
Chow$_\delta(X)$ be the set of all algebraic cycles of $X$
representing the homology class $\delta$.  We can embed $X$ into
some projective space ${\mathbb P}^n$, and then cycles of $X$ of
homology class $\delta$ become cycles of  ${\mathbb P}^n$ of
degree $d$ for some $d$. Hence Chow$_\delta(X)$ is embedded into
$\hbox{Chow}_d ({\mathbb P}^n)$, and  the variety Chow$_\delta(X)$
with the induced projective algebraic structure is called Chow
variety of algebraic cycles of $X$ of homology class $\delta$.
(See \cite{Ko} for further details.)

\subsection{Definition of Chow quotient}
Consider a reductive algebraic group action on a projective variety over the field of
complex numbers
$$G \times X \to X,$$
besides Mumford's GIT quotient, Kapranov-Sturmfels-Zelevensky (\cite{KSZ}) for toric varieties
and Kapranov (\cite{Ka}) in general, introduced the canonical Chow quotient.

\begin{defn} (Kapranov \cite{Ka}).
\label{defnChowQuotient} Let $x_0$ be a fixed {\it generic} point
of $X$ and $\delta$ be the induced homology class represented by
the cycle $\overline{G \cdot x_0}$. There is a Zariski open
subset, $U \subset X$, containing the point $x_0$, such that
$\overline{G \cdot x}$ represents the same homology class $\delta$
for all $x \in U$. Let Chow$_\delta(X)$ be the component of the
Chow variety of $X$ containing cycles of homology class $\delta$.
Then there is an embedding
$$\iota: U/G \to \hbox{Chow}_\delta (X)$$
$$ [G \cdot x] \to [\overline{G \cdot x}] \in \hbox{Chow}_\delta (X).$$
The Chow quotient, denoted by $X/\!/^{ch}G$, is defined to be 
the closure of $\iota (U/G)$.
\end{defn}

This definition is independent of the choice of the Zariski open
subset $U$. Note that the group $G$ acts on the Chow variety
Chow$_\delta(X)$ by moving the Chow cycles, and under this action,
the Chow quotient $\overline{\iota (U/G)}$ is contained in the
fixed point set.


\begin{rem}
\label{generic} By our assumption of the action, there is a
Zariski open subset $U_0$ such that points of $U_0$ are
isotropy-free. Let $U$ be the Zariski open subset in Definition
\ref{defnChowQuotient}.  Now, we point out that, in this paper,
unless otherwise indicated, by generic points we shall always mean
points in the Zariski open subsets $U_0 \cap U$.
\end{rem}

\subsection{The Chow family}
\label{chowfamily}

Let $$F \subset X \times (X/\!/^{ch}G) $$ be the family of
algebraic cycles over the Chow quotient $X/\!/^{ch}G$ defined by
$$F= \{ (x, Z) \in X \times (X/\!/^{ch}G) | x \in Z\}.$$
Then, we have a diagram
\begin{equation*}
\begin{CD}
F @>{\ev}>> X\\
@V{f}VV  \\
X/\!/^{ch}G
\end{CD}
\end{equation*}
where $\ev$ and $f$ are  the projections to the first and second factor,
respectively.  For any point $q \in X/\!/^{ch}G $, we will call the fiber,
$f^{-1}(q)$, the Chow fiber over the point $q$, and sometimes denote it by $F(q)$.

By \cite{Atiyah}, for a generic point $x \in X$, we have
$$\Phi_L (\overline{G \cdot x}) = \Phi_L (X).$$
That is, for a generic point $q  \in X/\!/^{ch}G$, we have
$$\Phi_L (\ev(F(q))) = \Phi_L (X).$$
In fact, this is  true  for any point $p  \in X/\!/^{ch}G$.

\begin{prop}
\label{sameimage}
 For any point $p \in X/\!/^{ch}G$, we have
$\Phi_L (\ev(F(p))) = \Phi_L (X).$
\end{prop}
\proof This is proved in Theorem (0.5.1) of \cite{BBS85}. We give
 slightly different arguments. By replacing $L$ by a large tensor
power, we may assume that all polytopes of the form $\Phi_L
(\overline{G \cdot x})$ (for some $x \in X$)  are  lattice
polytopes (see Mumford's Appendix to \cite{Ness}). Let $p \in
X/\!/^{ch}G$ be an arbitrary point and $q \in X/\!/^{ch}G$ be a
nearby generic point.
Since $p$ and $q$ are nearby points, we have that $\ev(F(p))$ and
$\ev(F(q))$ are nearby in the Hausdorff topology on the compact
subsets of $X$. Hence by the continuity of $\Phi_L$, the two {\it
compact lattice} polytopes $\Phi_L (\ev(F(p)))$ and $\Phi_L
(\ev(F(q)))$ are also nearby. Therefore, they must be equal.
\endproof

Each Chow fiber $F(q)$ as an algebraic cycle is of the form
$\sum_i m_i \overline{G \cdot x_i}$ where $m_i$ are nonnegative
integers and $G \cdot x_i$ are orbits of the top dimension. We
will call $\cup_i \overline{G \cdot x_i}$ the support of $F(q)$
and denote it by $|F(q)|$. Note that $\Phi_L (\ev(F(q))) = \Phi_L
(|F(q)|)$.

Observe from Proposition \ref{sameimage} that $\Phi_L(|F(q)|)$
gives a subdivision of the polytope $P$ by the subpolytopes
$\Phi_L (\overline{G \cdot x_i})$, any two of which either do not
meet or intersect along a proper face (see Theorem (0.5.1) of
\cite{BBS85}).

We will make the following technical assumption which says that
there are no two Chow fibers with exactly the same support but
different multiplicities. This assumption is needed for most of
the paper except \S 3.6 (Perturbing, translating, and
specializing) which is not affected.

\begin{assum}
\label{uniquesupport}  If $p \ne q$, then $|F(p)| \ne |F(q)|$.
\end{assum}

We do not know whether this assumption always holds, nor do we
know an example violating the assumption. When $X = {\mathbb P}^n$
and  $G$ is a subtorus of the dense open torus $(\C^*)^n \subset
{\mathbb P}^n$,
if a Chow fiber $F(q)= \sum_i m_i \overline{G \cdot x_i} $ is  a
minimal toric degeneration of the generic orbit closures (also
called extreme toric degeneration), then $m_i$ is {\it uniquely}
computed by the volume of certain simplex of $\dim G$ (Theorem
3.2, Chapter 8, \cite{GKZ}. See also \cite{KSZ}). Hence, it seems
that the assumption holds in this case.


It has been known that the Chow quotient dominates all GIT quotients (\cite{Ka}).
For any GIT quotient $X_C$ (\S \ref{GIT}), we shall use
$$\pi_C : X/\!/^{ch}G \rightarrow X_C$$
to denote the corresponding canonical projection.

\begin{rem}
\label{ChowtoGIT} In fact, let $U$ be any invariant open subset such that
the compact geometric quotient $U/G$ exists.  Then for any $q \in X/\!/^{ch}G$,
$F(q) \cap U$ is a single orbit in $U$ (Theorem 0.4, \cite{B-BS}). In particular,
let $\pi: X/\!/^{ch}G \rightarrow U/G$ be the projection, then
$\pi (q) = [F(q) \cap U] \in U/G$.
\end{rem}


Before we proceed further,
two conventions are needed.
\begin{enumerate}
\item The symbols
$[\overline{G \cdot x}]$ or $[\sum_i \overline{G \cdot x_i}]$ will be used to
denote a point in the Chow quotient $X/\!/^{ch}G$;
\item while $[G \cdot x]$ will be used to
denote a point in a given GIT quotient $X^{ss}/\!/G$.
\end{enumerate}

\subsection{Relation with the limit quotient}
 Recall that for a
fixed ample linearized line bundle $L$, the  moment map image $P =
\Phi_L (X)$ has a natural wall and chamber structure. The set of
GIT quotients associated to the moment map $\Phi_L$ is indexed by
the set ${\mathcal C}_L$ of all chambers. They form a finite
inverse system $\{X_C, f_{CD}| D \subset C  \in {\mathcal C}_L
\}$, that is, a finite set of projective varieties,  together with
canonical projective morphisms $f_{CD}$ among them.

\begin{defn} Let $\lim_{C \in {\mathcal C}_L } X_C$ be the inverse limit of
the system $\{X_C, f_{CD}| D \subset C  \in {\mathcal C}_L\}$.
($\lim_{C \in {\mathcal C}_L } X_C$ may be reducible, in general.)
The
 unique irreducible component of $\lim_{C \in {\mathcal C}_L} X_C$ that
contains the open subset $U/G$ is called the {\it limit quotient}
of $X$ by $G$ (associated with $L$), and is denoted by
$X/\!/^{lim}_L G$.
\end{defn}

It has been expected for quite a while that the Chow quotient
should in some way be related to the inverse limit of GIT
quotients.
Since it lacks a reference, we provide the necessary details
below.

\begin{lem}
\label{p=q} Suppose that Assumption \ref{uniquesupport} holds. Let
$p$ and $q$ be two points in $X/\!/^{ch}G$. Then $p=q$ if and only
if $\pi_C(p) = \pi_C (q)$ for all maximal (hence all) chambers
$C$.
\end{lem}
\proof Since all the projections, $\pi_C: X/\!/^{ch}G \rightarrow
X_C$, factor through GIT quotients defined by maximal chambers, we
only need to consider maximal chambers $C$. The necessary
direction is trivial. For sufficient direction, assume that $p\ne
q$. Then by Assumption \ref{uniquesupport}, $|F(p)| \ne |F(q)|$.
By Proposition \ref{sameimage},  we always have
$$\Phi_L (|F(p)|) =\Phi_L (|F(q)|) =\Phi_L (X).$$
Hence we can find $x \in |F(p)|$ and $y \in |F(q)|$ with
$$\dim G \cdot x = \dim G \cdot y = \dim G$$
such that
$$G \cdot x \ne G \cdot y$$ and
$$ \Phi(G \cdot x) \cap \Phi(G \cdot y) \ne \emptyset.$$
Now choose any maximal chamber $C$ such that
$$C \subset  \Phi(\overline{G \cdot x}) \cap \Phi(\overline{G \cdot y}).$$
(Here we have used the fact that $x$ is a regular point of the
moment map $\Phi$ if and only if $\dim G \cdot x  = \dim G$.)  It
follows that $x, y \in X^{ss}(C)$.

Then by Remark \ref{ChowtoGIT}, $\pi_C(p)= [G \cdot x]$ and
$\pi_C(q)= [G \cdot y]$. But since $C$ is a top chamber, all
points in $X^{ss}(C)$ are stable. Hence the two different orbits
$G \cdot x$ and  $G \cdot y$ will not be identified in the
geometric quotient  $X_C$.  Thus, we obtain
$$\pi_C(p)= [G \cdot x] \ne [G \cdot y] = \pi_C(q).$$
\endproof

\begin{thm}
\label{chow=limit} Suppose that Assumption \ref{uniquesupport}
holds. Then there is a natural birational projective
  morphism $\ell$ from $X/\!/^{ch}G $ to
$X/\!/^{lim}_LG$ which is also bijective. In particular,
$X/\!/^{ch}G $ is homeomorphic to
$X/\!/^{lim}_LG$. 
\end{thm}
\proof First of all, since $X/\!/^{ch}G $ maps naturally to all
GIT quotients $X_C$, we have that $X/\!/^{ch}G $ maps naturally to
$\lim_{C \in {\mathcal C}} X_C$, and the map is an isomorphism
when restricted to $U/G$. Hence the image is contained in
$X/\!/^{lim}_LG$. This gives a birational surjective projective
morphism
$$\ell: X/\!/^{ch}G \rightarrow X/\!/^{lim}_L G $$
It suffices to show that $\ell$ is injective. For any two points
$p, q \in X/\!/^{ch}G$, if $\ell (p) = \ell (q)$, then $\pi_C (p)
= \pi_C (q)$ for all $C \in {\mathcal C}$. By Lemma \ref{p=q}, we
have $p=q$.

$\ell$ is homeomorphic because a bijective continuous map between
two compact Hausdorff spaces is a homeomorphism.
\endproof


Consequently, we have that $X/\!/^{lim}_L G$ is independent of the
ample line bundle $L$, at least up to homeomorphism. This may look
surprising. But when $G$ is a torus -- and this is what we deal
with in this paper, changing the line bundle $L$ will deform the
moment map image; under this deformation, a very few GIT quotients
change slightly: a few older GIT quotient disappear, while a few
new ones emerge, but the limit quotient remains homeomorphic to
each other. The same can not be said when $G$ is a general
reductive group. In fact, in the general case, even if we take
limit over all ample line bundles, we suspect that Theorem
\ref{chow=limit} should be false.

\subsection{Ample line bundles over the Chow quotient}

Fix a very ample linearized line bundle $L$. Replacing $L$ by a
large tensor power, if necessary, we may assume that $L$ descends
to a very ample line bundle $L_C$ over the GIT quotient $X_C$,
simultaneously for all chambers $C \in \cC_L$.

Define a line bundle $L_{ch}$ over $X/\!/^{ch}G$ by setting
$$L_{ch} = \otimes_{\hbox{maximal} \; C \in \cC_L} \pi_C^* L_C.$$
We will show that this is an ample line bundle over $X/\!/^{ch}G$.

\begin{lem}
\label{curvelemma}
For any curve $Z \subset X/\!/^{ch}G$, there exists a maximal chamber $C$ such that
$\pi_C (Z)$ is a nontrivial curve in $X_C$.
\end{lem}
\proof Pick two distinct points $q_1$ and $q_2$ in $Z \subset
X/\!/^{ch}G$. If $\pi_C (q_1) = \pi_C (q_2)$ for all chambers $C$,
then $q_1 = q_2$ by Lemma \ref{p=q}. Hence there is a chamber $C$
such that $\pi_C (q_1) \ne \pi_C (q_2)$.
 This shows that $\pi_C (Z)$ is a nontrivial curve in the GIT quotient
$X_C$.
\endproof

\begin{prop}
$L_{ch}$ is ample over $X/\!/^{ch}G$.
\end{prop}
\proof
Take any curve $Z$ in $X/\!/^{ch}G$. By the previous lemma, there is a maximal chamber
$C_0$ such that $\pi_{C_0} (Z)$ is a nontrivial curve in $X_{C_0}$. We have
$$L_{ch} \cdot [Z] = \sum_{\hbox{maximal} \; C \in \cC_L} \pi_C^* L_C \cdot [Z]$$
$$= \sum_{\hbox{maximal} \; C \in \cC_L}  L_C \cdot [\pi_C (Z)].$$
All $L_C \cdot [\pi_C (Z)] \ge 0$, and $L_{C_0} \cdot [\pi_{C_0}
(Z)] > 0.$ Hence $L_{ch} \cdot [Z] >0$. Now observe that $L_{ch}$
is semi-ample by the definition, hence the above positivity
implies that it can not contract any curve. Therefore, it is
ample.
\endproof

\begin{rem}
If we define $L_{ch}'$ by setting
$$L_{ch}' = \sum_{\hbox{all } \; C \in \cC_L} \pi_C^* L_C,$$
then, this is also an ample line bundle. The proof is the same.
\end{rem}

\begin{rem}
Choose a $K$-invariant symplectic form $\omega$ in $[c_1(L)]$. For
each top chamber $C$, choose an interior point ${\bf r}_C \in C$.
Then via symplectic reduction, we obtain a symplectic form
$\omega_{{\bf r}_C}$ on the quotient $X_C$. Let $\pi_C^*
\omega_{{\bf r}_C}$ be the pullback of the form on $X/\!/^{ch}G$.
Then, add them up for all $C$, we obtain
$$\sum_C \pi_C^* \omega_{{\bf r}_C}$$
which is a symplectic form on (the smooth locus of) $X/\!/^{ch}G$
because the map $X/\!/^{ch}G \rightarrow X/\!/^{lim}G$ does not
contract anything.
\end{rem}

\subsection{Perturbing, translating and specializing}

From this subsection till \S 3.8, we do not suppose that
Assumption \ref{uniquesupport} holds.

A point of $X$ is said to be isotropy-free if its isotropy subgroup is the identity subgroup.

\begin{thm}
\label{chowfiberAG}
Let $x$ and $y$ be two points in $X$ such that
$\dim G \cdot x = \dim G \cdot y = \dim G$.
Then the points $x$ and $y$  belong to the same
Chow fiber $F(q)$ for some $q \in X/\!/^{ch}G$ if and only if
there is a generic holomorphic map from the complex unit disk $\Delta= \{z| |z| <1 \}$
to $X$
$$\varphi: \Delta \rightarrow X$$
with  $\varphi (0) = x$ and
a holomorphic map from the punctured disk $\Delta^*=\Delta \setminus
\{0\}$ to $G$
 $$g: \Delta^* \rightarrow G$$
such that $$y=\lim_{t \to 0} g(t) \cdot \varphi (t).$$
\end{thm}
\proof (For an illustration of this theorem, see Figure 2.)

To prove the sufficient part, we can assume that, for $t \in
\Delta^*=\Delta \setminus \{0\}$, the orbits $G \cdot \varphi(t)$
are generic ( see Remark \ref{generic}),
 hence we obtain
a well-defined  holomorphic map
$$\tilde\varphi: \Delta^* \rightarrow X/\!/^{ch}G$$ with
$$\tilde\varphi (t) = [\overline{G \cdot \varphi(t)}].$$
Now,  because $\varphi (0) = x$, we  have
$$\lim_{t \to 0} \tilde\varphi  (t) = F_x,$$
where $F_x$ is a Chow fiber that contains $G \cdot x$.
On the other hand, we have that
$$\tilde\varphi  (t) = [\overline{G \cdot  \varphi(t)}] = [\overline{G \cdot g(t)\cdot  \varphi(t)}].$$
Now take the limit $t \to 0$, since $y = \lim_{t \to 0} g(t) \cdot \varphi (t),$
we obtain that
$$\lim_{t \to 0} \tilde\varphi  (t) = F_y,$$
where $F_y$ is a Chow fiber that contains $G \cdot y$.
Hence $F_x = F_y$.

Conversely, assume that the two points $x$ and $y$  belong to the same
Chow fiber, $F(q)$, for some $q \in  X/\!/^{ch}G$.

First, since $U_0 \cap U$ is Zariski open (see Remark
\ref{generic}), using Luna's \'etale slice around the point $x$
(see \cite{Luna}), we can choose a  holomorphic map
$$\varphi: \Delta \rightarrow  X$$
with $\varphi(0) =x$ and $\varphi (t) \in U_0 \cap U$ for all $t
\ne 0$ such that $[\overline{G \cdot \varphi (t)}] \in
X/\!/^{ch}G$  approaches the Chow point $q$ as $t$ goes to $0$.
Note in particular that the orbit $G \cdot \varphi(t)$ is an
isotropy-free generic point for  $t \ne 0$.



Next, for the point $y$, there is an invariant open subset $U_y$
such that $U_y$ contains $y$ and the compact geometric GIT
quotient $U_y/G$ exists. For example, take any top chamber $C$
that is contained in the polytope $\Phi_L (\overline{G \cdot y})$
(which is also top dimensional), and let $U_y = X^{ss}(C)$.


Now, for $t \ne 0$, consider the orbits $$[G \cdot \varphi(t)]_y \in U_y/G$$
as points in the geometric quotient $U_y/G$, and $$[\overline{G \cdot \varphi(t)}]\in X/\!/^{ch}G$$
as points in the Chow quotient.
We have the following  diagram
\begin{equation*}
\begin{CD}
[\overline{G \cdot \varphi(t)}] @>{t \to 0}>> q \\
@V{\pi}VV  @V{\pi}VV \\
[G \cdot \varphi(t)]_y   @>{t \to 0}>> [F(q) \cap U_y]_y
\end{CD}
\end{equation*}
where the horizonal arrows are taking limits and the down arrows are the projection
morphism from Chow quotient $X/\!/^{ch}G$ to the GIT quotient $U_y/G$.
The top horizontal arrow $\lim_{t \to 0} [\overline{G \cdot \varphi(t)}] = q$ holds because of the choice of $\varphi$.
The down arrows hold because of Remark \ref{ChowtoGIT}. Now by the continuity of $\pi$, we obtain
the bottom horizontal arrow
$$\lim_{t \to 0} [G \cdot \varphi(t)]_y  =[F(q) \cap U_y]_y.$$

But $G \cdot y \subset F(q) \cap U_y$. Hence $G \cdot y = F(q) \cap U_y$ (Remark \ref{ChowtoGIT}).
Therefore we obtain
$$\lim_{t \to 0} [G \cdot \varphi(t)]_y =[G \cdot y]_y \in U_y/G. \;\;\;\;\;\;\;\;\;\;\; (*)$$

For the point $y$, there is an analytic slice by Luna's \'etale
slice Theorem (\cite{Luna}), $y \in S_y \subset U_y$, such that
$S_y$ meets transversally, at a unique point, every $G$-orbit in
the open subset $W_y = G \cdot S_y$. Hence there is $\delta >0$
such that
$$\{G \cdot \varphi (t) \cap S_y | 0< |t| < \delta \}$$
is a holomorphic curve in $S_y$. Let
$G \cdot \varphi (t) \cap S_y = g(t) \cdot \varphi (t) \in S_y$ for $0< |t| < \delta$.
Because $\varphi (t)$ is isotropy-free and the orbit $G \cdot \varphi (t)$ meets the slice $S_y$ transversally
at a unique point,
$g(t)$ is holomorphically uniquely determined for each
$ 0< |t| < \delta$. In other words, $g(t)$ is equivalent to the holomorphic map
$$\{t | 0< |t| < \delta\} \rightarrow S_y$$
$$ t \rightarrow g(t) \cdot \varphi (t).$$
Hence we have that
$$g: \{z | 0< |z| < \delta\} \rightarrow G$$
is a holomorphic map.
Let $y' \in S_y$ be the limit of $g(t) \cdot \varphi(t)$ as $t$ approaches 0.
Since $[G \cdot \varphi (t)]_y = [G \cdot g(t) \cdot \varphi (t)]_y$,
we obtain, by the identity (*), that
$$[G \cdot y]_y = \lim_{t \to 0} [G \cdot g(t) \cdot \varphi(t)]_y = [G \cdot y']_y.$$
This implies that $y'=y$ because $y, y' \in G \cdot y \cap S_y$.
Now,  by a suitable parameter change, we may assume $\delta =1$.

This completes the proof.
\endproof

\begin{rem}
In the proof, we require that $\varphi (t)$ has the trivial
isotropy subgroups. This is not necessary. But it simplifies proofs and also
some applications as we will see in  \S \ref{section:toppush}.
\end{rem}

\begin{rem}
\label{topnature}
This theorem is {\it topological} in nature.
In other words, if we replace all the word ``holomorphic (maps)''
by ``continuous (maps'', the theorem and its proof remain unchanged.
The theorem also holds when we replace ``Chow'' by ``Hilbert''.
This, again, is due to (shows)
the topological nature of the theorem.
\end{rem}

\begin{rem} Now, a few words on the terminologies. The effect of the map $\varphi: \Delta \to X$
is to push the point $x = \varphi(0)$ to  general positions, $\varphi (t)$, $t \ne 0$; then
the group elements $g(t)$ translate the points $\varphi (t)$ along the group orbits
to  positions close to $y$, allowing the final desired specialization.
This motivates us to use the descriptive terms
``perturbing,  translating and specializing''. We sometimes abbreviate it as ``P.T.S''.
\end{rem}

\begin{defn}
We say that a point $x$ of $X$ can be perturbed (to general positions),
translated (along $G$-orbits),
 and specialized to a point $y$ of $X$ if they have the relation as
 described in Theorem \ref{chowfiberAG}.
In this case, we will write $x \to_G y$.
\end{defn}

\begin{rem}
\label{symetric}
Let $X_{(0)}$ be the set of all points of $X$ whose isotropy subgroups are finite.
Then P.T.S. formulation defines a relation on $X_{(0)}$.
By the above theorem,  this relation is equivalent to the relation defined by:
$x$ and $y$ ``belong to the same Chow fiber'', which is obviously symmetric.
This can also be seen directly. Assume that we have $x \to_G y$, that is,
 there is a generic holomorphic map from the unit disk to $X$
$$\varphi: \Delta \rightarrow X$$
with $\varphi (0) = x$, and a holomorphic map from the punctured disk
$ \Delta^*= \Delta  \setminus \{0\}$ to $G$
 $$g: \Delta^* \rightarrow G,$$
such that $$y=\lim_{t \to 0} g(t) \cdot \varphi (t).$$
Define $\psi (t) = g(t) \varphi (t)$. Then, this can be extended to
a holomorphic map from the unit disk to $X$
with $\psi (0) = y$. Let $h(t) = g^{-1}(t)$. Then this defines a holomorphic
map from the punctured disk $\Delta^*=\Delta \setminus \{0\}$ to $G$.
Clearly, we have
$$x=\lim_{t \to 0} h(t) \cdot \psi (t).$$
That is, $y \to_G x$
\end{rem}

However,  an orbit of points of  $X_{(0)}$ may belong to several different Chow fibers.
(This may happen only when the orbit is not generic).
Hence  P.T.S. fails to be transitive, and thus
is not an equivalent relation on $X_{(0)}$.

\subsection{Point configurations on ${\Bbb P}^n$}

In this section, we draw some special consequences for the diagonal action of
$\PGL(n+1, {\Bbb C})$ on $({\Bbb P}^n)^m$.

Theorem \ref{chowfiberAG} holds for all torus action, and in
particular, holds for the action of the maximal torus $H = ({\Bbb
C}^*)^m/\Delta$ on $\Gr(n+1,{\Bbb C}^m)$, the Grassmannian  of
$n+1$-dimensional planes in ${\Bbb C}^m$  ($n+1 <m$),  where
$\Delta$ is the diagonal subgroup of $({\Bbb C}^*)^m$ which acts
trivially on $\Gr(n+1,{\Bbb C}^m)$. Using the Gelfand-MacPherson
correspondence, we can transform the properties of the  $H$-action
on $\Gr(n+1,{\Bbb C}^m)$ to the corresponding properties of the
$G$-action on $({\Bbb P}^n)^m$ where $G=\PGL(n+1)$.

The GM correspondence can be seen as follows. We represent a point
of $\Gr(n+1,{\Bbb C}^m)$ as a matrix $A$ of size $(n+1) \times m$.
Write $A$ in column vectors, $A = (a_1, \ldots, a_m)$, with $a_i
\in \C^{n+1}$. $A$ is of full rank and we assume that all $a_i \ne
0$. Let $U$ be the set of all such matrices.  The group $\GL(n+1)$
acts on $A \in U$ from the left. The group $(\C^*)^m$ acts on $A$
from the right by multiplying the column vectors componentwise.
Take the quotient  of $U$ by $\GL(n+1)$ first, we obtain the
Grassmannian $\Gr(n+1,{\Bbb C}^m)$ with the residual $H = ({\Bbb
C}^*)^m/\Delta$- action. Take the quotient of $U$ by $(\C^*)^m$
first, we obtain $({\Bbb P}^n)^m$ with the residual diagonal
action of $G=\PGL(n+1)$.  This establish a correspondence between
$H$-orbits on $\Gr(n+1,{\Bbb C}^m)$ and $G$-orbits on $({\Bbb
P}^n)^m$: they all correspond to $\GL(n+1) \times_{\C^*}
(\C^*)^m$-orbits on $U$. Then, taking quotient in stages, we get a
natural correspondence between GIT quotients as well as Chow
quotients (see \cite{Ka} for details. See also \cite{Hu04} for a
generalization of GM correspondence to the product of
Grassmannians).

Thus, by the GM correspondence and in terms of point
configurations on $X=({\Bbb P}^n)^m$, Theorem \ref{chowfiberAG}
reads

\begin{thm}
\label{pushTrans}
Let $\underline{x}$ and $\underline{y}$ be two points in $X=({\Bbb P}^n)^m$ such that
$\dim G \cdot \underline{x} = \dim G \cdot \underline{y} = \dim G$.
Then, the points $\underline{x}$ and $\underline{y}$  belong to the same
Chow fiber, $F(q)$, for some $q \in X/\!/^{ch}G$ if and only if
there is a generic holomorphic map from the complex unit disk $\Delta= \{z| |z| <1 \}$
to $X$
$$\varphi: \Delta \rightarrow X$$
with $\varphi (0) =\underline{x}$ and
a holomorphic map from the punctured disk $\Delta^*=\Delta \setminus
\{0\}$ to $G$
 $$g: \Delta^* \rightarrow G$$
such that $$\lim_{t \to 0} g(t) \cdot \varphi (t) = \underline{y}.$$
\end{thm}



The theorem bears an interesting corollary in the  case of $({\Bbb P}^1)^m$.
Let $\{1, \ldots,m\} = J \cup J^c$,  where $J$ is a subset containing at least two elements and $J^c$ is
the complement.

\begin{thm}
\label{pointConfiginP1} Let $\underline{x}=(x_1, \ldots, x_m)$ be
a point in $({\Bbb P}^1)^m$ such that its isotropy  subgroup is
trivial, that is, at least three points are distinct. Assume that
$\underline{x}_J$, the points with indexes in $J$,
 coincide, but $\underline{x}_{J^c}$ are all distinct.
Let $\underline{y}= \lim_{t \to 0} g(t) \cdot \varphi (t)$ such that $\underline{y}_J$ are all distinct.
Then $\underline{y}_{J^c}$ must coincide.
\end{thm}

\proof By Theorem \ref{chowfiberAG}, the algebraic cycle
$$[\overline{G \cdot \underline{x}}] + [\overline{G \cdot \underline{y}}]$$
lies in a Chow fiber $F(q)$. Via Kapranov's isomorphism
(\cite{Ka}) between $({\Bbb P}^1)^m/\!/^{ch}G$ and
$\overline{M}_{0,n}$ (the moduli space of stable $n$-pointed
rational curves),   $F(q)$ corresponds to a stable $n$-pointed
rational curve $C$. The curve $C$ contains two components, one
component, corresponding to the  cycle $[\overline{G \cdot
\underline{x}}]$, contains the distinct points
$\underline{x}_{J^c}$; and another component, corresponding to the
cycle $[\overline{G \cdot \underline{y}}]$, contains the distinct
point $\underline{y}_J$.  Since $\underline{x}_{J^c}$ and
$\underline{y}_J$ together give $n$ distinct points on the stable
$n$-pointed rational curve $C$, we see that there should be no
other components in $C$. Therefore $\underline{y}_{J^c}$ must
coincide, and the two components are glued by joining the point
$\underline{x}_J$ with the point $\underline{y}_{J^c}$.
\endproof

This elementary and interesting ``new'' phenomenon seems elusive
before the discovery of Theorem \ref{chowfiberAG}. It can, indeed,
be explained by the elementary complex analysis. We think of
${\Bbb P}^1$ as the extended plane. As earlier, using $ \varphi
(t)$,  we can separate $\underline{x}_J$ to distinct points
$\underline{x}_J (t)$ ($t \ne 0$) around $x_j$ ($j \in J$). For
simplicity,  we may assume that $x_j = 0$. Now we want to apply a
one-parameter curve $g (t)$ in $\PGL(2)$ such that in the limit, $
g(t) \cdot \underline{x}_J (t)$ get separated. Recall that linear
transformations are made of translations, rotations, homotheties,
and inversions. Among the four kinds, only inversion will do the
work. Hence we may assume that $g(t)$ are inversions. Take a small
neighborhood $D$ of  $x_j=0$, inversions $g(t)$ will amplify $D$
to a large neighborhood at infinity, and simultaneously shrink the
complement of $D$ into a small neighborhood around 0. Hence after
taking limit, the neighborhood $D$ expands to the whole extended
plane, and in the mean time, the complement of $D$ collapses to
the single point 0 --- and this explains why the points $x_{J^c}$
collide into a single point in the end.

\begin{exmp}
As a concrete example, take a set of points of  $({\Bbb P}^1)^m$
represented by a $2 \times m$ matrix
$$
\left(
\begin{array}{ccccccc}
            a & a &  \cdots & a & a_{j+1} &  \cdots & a_m \\
            b & b &  \cdots & b & b_{j+1} &  \cdots & b_m
\end{array}
\right) \in ({\Bbb P}^1)^m
$$
with $b \ne 0$ such that the first $j$ points coincide, and the rest
$$\left(
\begin{array}{ccccccc}
             a_{j+1} &  \cdots & a_m \\
             b_{j+1} &  \cdots & b_m
\end{array}
 \right)
$$
is sufficiently general.
Perturb $$
\left(
\begin{array}{ccccccc}
            a & a &  \cdots & a & a_{j+1} &  \cdots & a_m \\
            b & b &  \cdots & b & b_{j+1} &  \cdots & b_m
\end{array}
\right)
$$
to a general position $\varphi (t)$ as
$$
\left(
\begin{array}{ccccccc}
            e^t a & e^{2t}a &  \cdots &  e^{jt}a & a_{j+1} &  \cdots & a_m \\
            b & b &  \cdots & b & b_{j+1} &  \cdots & b_m
\end{array}
\right).
$$
Let $g(t)$ be given by
$$\left(
\begin{array}{ccccccc}
             {1 \over t} &  - {a \over b} {1 \over t} \\
             0 &  1
\end{array}
 \right).
$$
Then $g(t) \cdot \varphi (t)$ is
$$
\left(
\begin{array}{ccccccc}
            {e^t-1 \over t} a & {e^{2t}-1 \over t}a &  \cdots & {e^{jt}-1 \over t}a &
{1 \over t} (a_{j+1} -{a \over b} b_{j+1})&  \cdots & {1 \over t}(a_m -{a \over b} b_m) \\
            b & b &  \cdots & b & b_{j+1} &  \cdots & b_m
\end{array}
\right).
$$
Let $t$ go to 0, we obtain a new set of points
$$
\left(
\begin{array}{ccccccc}
            a & 2a &  \cdots & ja & 1 &  \cdots & 1 \\
            b & b &  \cdots & b & 0 &  \cdots & 0
\end{array}
\right)
$$
where,  as predicted in the theorem,  the first $j$ points get separated, and the rest collide at a single point.
\end{exmp}

\begin{exmp} We can also think of
$$
\left(
\begin{array}{ccccccc}
            a & a &  \cdots & a & a_{j+1} &  \cdots & a_m \\
            b & b &  \cdots & b & b_{j+1} &  \cdots & b_m
\end{array}
\right) \in \Gr(2, m)
$$
as a point in the Grassmannian $ \Gr(2, m)$. We can have the same perturbation
$$
\varphi (t) = \left(
\begin{array}{ccccccc}
            e^t a & e^{2t}a &  \cdots &  e^{jt}a & a_{j+1} &  \cdots & a_m \\
            b & b &  \cdots & b & b_{j+1} &  \cdots & b_m
\end{array}
\right).
$$
But $\varphi (t)$,  as a point in $ \Gr(2, m)$, is the same as $g(t) \cdot \varphi (t)$,
which, in terms of matrix,  is
$$
\left(
\begin{array}{ccccccc}
            {e^t-1 \over t} a & {e^{2t}-1 \over t}a &  \cdots & {e^{jt}-1 \over t}a &
{1 \over t} (a_{j+1} -{a \over b} b_{j+1})&  \cdots & {1 \over t}(a_m -{a \over b} b_m) \\
            b & b &  \cdots & b & b_{j+1} &  \cdots & b_m
\end{array}
\right).
$$
Now, let $h(t) \in ({\Bbb C}^*)^m$ be given by
$$
(1, \ldots, 1, t {1 \over a_{j+1} -{a \over b} b_{j+1}}, \ldots, t {1 \over a_m -{a \over b} b_m}).
$$
Then $h (t) \cdot  \varphi (t)$ becomes
$$
\left(
\begin{array}{ccccccc}
            {e^t-1 \over t} a & {e^{2t}-1 \over t}a &  \cdots & {e^{jt}-1 \over t}a &
1 &  \cdots & 1 \\
            b & b &  \cdots & b & t {1 \over a_{j+1} -{a \over b} b_{j+1}} b_{j+1} &  \cdots & t {1 \over a_m -{a \over b} b_m} b_m
\end{array}
\right).
$$
Let $t$ go to zero, we again obtain
$$
\left(
\begin{array}{ccccccc}
            a & 2a &  \cdots & ja & 1 &  \cdots & 1 \\
            b & b &  \cdots & b & 0 &  \cdots & 0
\end{array}
\right).
$$
\end{exmp}
As one can see, although equivalent, it is sometimes easier to describe
 the properties of a matrix as a configuration of points
on ${\Bbb P}^1$ than as a 2-plane in ${\Bbb C}^m$.

\subsection{A geometric interpretation of Chow quotients of higher Grassmannians}

As far as moduli spaces of point configurations on ${\Bbb P}^n$ are concerned, the case of $n=1$
is very special. Here, by higher Grassmannians, we mean $\Gr(n,m)$ with $n >2$, $m >5$. Again,
the Chow quotients of higher Grassmannians correspond to Chow quotients of $({\Bbb P}^n)^m$ with $n >1, m >5$.

Take a point configuration, $\underline{x}=(x_1, \ldots, x_m)$, in ${\Bbb P}^n$,
such that its automorphism group is trivial. By Theorem \ref{pushTrans}, up to projective transformations,
there are only finitely many points, including $\underline{x}$ itself,
 $$\underline{x}_1, \ldots, \underline{x}_l,$$
which can be obtained from $\underline{x}$ by perturbing,
translating and specializing.

\begin{defn}
\label{defn:geom}
We will call $\underline{x}_1, \ldots, \underline{x}_l$, or the union of their $G$-orbits,
a stable configuration of $m$-points in ${\Bbb P}^n$.
\end{defn}

Then by Theorem \ref{pushTrans},  we have

\begin{thm}
\label{thm:geom}
The Chow quotient of $({\Bbb P}^n)^m$ by the group $\PGL(n+1, {\Bbb C})$ parameterizes
stable configurations of $m$-points in ${\Bbb P}^n$.
\end{thm}

\begin{rem}
Note that from the original definition, the Chow quotient is defined by taking the closure of $\iota (U/G)$.
Taking closure usually does not provide further information on the boundary points.
Definition \ref{defn:geom} and Theorem \ref{thm:geom}, relying on the computable
{\it perturbing-translating-specializing} formulation, fill the gap to certain degree.
\end{rem}

Ideally, it would be desirable to specify how the orbit closures,
$\overline{G \cdot \underline{x}_i}$, are glued together like the
case of $n=1$ (cf. Theorem \ref{pointConfiginP1} and the proof).
Despite the intimidating combinatorial complexity, we hope that
the ideas presented here will lead to, in a forthcoming paper,
 a much better understanding
of the Chow quotients of higher Grassmannians\footnote{After
posting this paper on ArXiv, we saw Paul Hacking's paper
\cite{Hacking} where he, geometrically interprets the Chow
quotient as the moduli space of stable log pairs. Later, Keel and
Tevelev posted  \cite{KT} which contains some results similar to
Hacking's.}.

\bigskip
\section{Chow Quotient:  symplectic approach}

Throughout the rest of the paper, we suppose that Assumption
\ref{uniquesupport} holds.

We begin with discussions on two questions (\S\S 4.1, 4.2) that
motivate the topic of this section.

\subsection{Symplectic reductions for the Chow quotient?}
We have known that a GIT quotient can be identified with various symplectic reductions.
Put it differently, a GIT quotient carries many (stratified) symplectic structures.
The connection is established by the theory of moment map. In the same vein, we may  ask:
{\it what are ``symplectic reductions'' for the Chow quotient?}
This is a natural question. It admits  inspiring approximations in
the following interesting cases.

Consider the diagonal action of $\PGL (2)$ on $({\bf P}^1)^{n+3}$.
Kapranov proved that the Chow quotient of this action is
$\overline{M}_{0,n}$, the moduli space of stable $n$-pointed
rational curves. In \cite{Hu99a}, we give a family of
``symplectic'' constructions of $\overline{M}_{0,n}$ using {\it
stable polygons with prescribed side lengths}. To say it
differently,  moduli spaces of stable polygons are in some sense
``symplectic reductions'' for the Chow quotient
$\overline{M}_{0,n}$.

This case is rather special, in that (stable) polygons play indispensable roles.
But, for the general case, it does inspire us to introduce the following new notion,
{\it stable orbits with fixed momentum
charges}, to take the role of stable polygons with fixed side lengths.
To further motivate  the precise description of these new objects, we next explore intuitively
what they should mean geometrically.

\begin{rem}
In what follows, the word ``orbit'' will always refer to ``$K$-orbit''.
When another group is involved, we will specify the group,
e.g., $G$-orbits.
\end{rem}

\subsection{Geometrically meaningful compactification}
\label{gmc}

Take a compact form $K$ of $G$. Let ${\bfmit k}$ be the Lie algebra of $K$.
As mentioned earlier,
we will focus on torus actions only.
Let $$\Phi: X \to {\bfmit k}^*$$
be a moment map for the $K$-action on $X$.
Pick any point  ${\bf r} \in \Phi(X)$.
Orbits in $\Phi^{-1} ({\bf r})$ will be said to have the momentum charge ${\bf r}$.
To keep with the theme of the rest of the  paper,
 we will, from now and on, frequently call symplectic reduction,
$$X_{\bf r}=\Phi^{-1} ({\bf r})/K,$$ the moduli space
of orbits with momentum charge ${\bf r}$. Let $X_{\bf r}^\circ
\subset X_{\bf r}$ be the moduli space of {\it generic orbits}
with momentum charge ${\bf r}$. Here an orbit $O = K \cdot x$ is
said to be generic if $\Phi(\overline{G \cdot x})$ equals the
whole polytope $\Phi(X)$. For example, orbits through the points
in the open subset $U$ of Definition \ref{defnChowQuotient} are
generic. Thus, $X_{\bf r}^\circ$ contains an open subset that is
homeomorphic to $U/G$, and is itself an open variety in $X_{\bf
r}$. From the definition, the orbits in the complement $X_{\bf r}
\setminus X_{\bf r}^\circ$,  measured by moment map image, are of
smaller size than those of generic ones. In the spirit of
geometric moduli problem, it is natural to ask for
compactifications of $X_{\bf r}^\circ$ with the following two
desirable characteristics:
\begin{enumerate}
\item the added boundary points should have natural geometric meanings and;
\item  the limiting geometric objects should be of the same size as the generic ones.
\end{enumerate}

To this end, we have proposed to add ``stable orbits'' as boundary points.
 So,  what are {\sl stable $K$-orbits}? First,
$K$-orbits through generic points in $\Phi^{-1} ({\bf r})$
are automatically considered to be {\it stable}. When a family of
generic orbits degenerate to a special orbit in $X_{\bf r}$, we
can imagine it as some kind of collision occurs, resulting orbits
of smaller sizes. In the case of
 spatial polygons, this means some edges become positively parallel (pointing to the same
direction). To get stable polygons, we introduce ``bubble''
polygons {\it with certain fixed side lengths} (\cite{Hu99a}). In
our current situation, what is needed is to introduce ``bubble''
orbits {\it with certain fixed momentum charges}. To know what
momentum charges to work with, some choice are to be made, just
like in the case of stable polygons, where we have to make choices
of side lengths. The detail is to be explicitly spelled out in the
subsequent section.

\begin{rem}
The moduli space ${\bfmit M}_\gamma$ of stable $K$-orbits to be
constructed below answers the question of \S \ref{gmc} quite
successfully. It is the interpretation of ${\bfmit M}_\gamma$ as a
"symplectic reduction" for the Chow quotient still unsatisfactory,
although it is obviously related. But, to find a "symplectic
reduction" for the Chow quotient is one of the motivations that
get this project started.
\end{rem}

\subsection{ Momentum charges}

Recall that we have the Chow family  $$F \subset X \times (X/\!/^{ch}G) $$
with the following diagram
\begin{equation*}
\begin{CD}
F @>{\ev}>> X\\
@V{f}VV  \\
X/\!/^{ch}G
\end{CD}
\end{equation*}
where $\ev$ and $f$ are  the first and second projection,
respectively.

For each $q \in X/\!/^{ch}G$, the support of the Chow fiber $F(q)$
is a union, $\cup_i \overline{G \cdot x_i}$, of orbit closures
with $\dim G \cdot x_i = \dim G$ for all $i$. The moment map image
of each orbit closure in $F(q)$ is a subpolytope of $\Phi(X)$, and
by the virtue of Proposition \ref{sameimage}, they all together
form a subdivision of $\Phi(F(q)) =  \Phi_L (X)= P_L$.

\begin{defn}
\label{subdivision}
 A coherent subdivision of $P_L$ is the
collection of top dimensional subpolytopes $\Phi_L (\overline{G
\cdot x_i})$ where $\cup_i \overline{G \cdot x_i} = |F(q)|$ for
some $q \in X/\!/^{ch}G$.
\end{defn}

  There are only finitely many
such polytopal subdivisions. We will use letter $S$ to denote such
a subdivision. We point out here that in this paper we will always
consider the subdivision $S$ as the {\it collection of the top
dimensional subpolytopes that occur in the subdivision.} Also, in
this paper, only coherent subdivisions will be considered, so we
will drop the word "coherent".

\begin{defn}
The set of all subdivisions of the form,
$$\Phi(F(q)), q  \in X/\!/^{ch}G,$$
will be denoted by ${\mathcal S}$. There is a partial order on the
set ${\mathcal S}$. For any two elements $S, S' \in {\mathcal S}$,
we say that $S < S'$ if $S$ is refined by $S'$. Under this partial
order, the poset ${\mathcal S}$ has a unique minimal element,
namely the (non-subdivided) polytope $P = \Phi(X)$.
\end{defn}

Fixed a general point ${\bf r}$ in $\Phi (X)$. For every polytopal
subdivision $S \in {\mathcal S}$, we choose a set of points,
$$\{ {\bf r}_D \in \; \hbox{the interior of } \; D | D \in S.\}$$
(Recall here that $D$ is of top dimension.) In other words, we
have an injective function,
$$\gamma_S: S \rightarrow \Phi (X),$$
 from the set of subpolytopes of $S$ to $\Phi(X)$,
 by sending
a polytope $D \in S$ to a point ${\bf r}_D $ in the interior
$D^\circ$ of $D$. Let $\gamma$  denote the collection of all the
above choices $\{\gamma_S: S \in {\mathcal S}\}$.

\begin{defn}
\label{charges} $\gamma$ is called an admissible set of momentum
charges with the principal charge ${\bf r}$ if the following
conditions are satisfied.
\begin{enumerate}
\item (principal main charge.) $\gamma_P (P) = {\bf r}$; \item
(local main charge.) Let a subdivision $S \in {\mathcal S}$ refine
another subdivision $S' \in {\mathcal S}$. Let $D \in S$ be
contained in $D' \in S'$, and $D$ contains $\gamma_{S'} (D')$.
Then $\gamma_{S} (D) = \gamma_{S'} (D')$. In particular, for any
subdivision $S$, we must have $\gamma_S  (D_{\bf r}) = {\bf r}$,
where $D_{\bf r}$  is the unique subpolytope in $S$ that contains
the original charge ${\bf r}$; \item  (compatibility.)
  For any two subdivisions $S$ and $S'$, if $D$ appears in both $S$ and $S'$, then
$$\gamma_S(D) = \gamma_{S'}(D).$$
\end{enumerate}
\end{defn}

\begin{rem}
Notice that (2) implies that if the subdivision $S \in {\mathcal
S}$ is refined by another subdivision $S' \in {\mathcal S}$, then
${\rm Image} (\gamma_{S}) \subset {\rm Image} (\gamma_{S'}).$ Note
also that in the sense of  (2), every given polytope $D$ admits a
(local) main charge $\gamma_S (D)$, and in particular, ${\bf
r}=\gamma_P (P)$ is the  global main charge.
\end{rem}

In this paper, only admissible set of
momentum charges will be considered. So, we will drop the word ``admissible''.
It worths to point out that
$\gamma$  is analogous to the choices $({\bf r}, \{{\bf \epsilon}\})$ of side lengths
in the case of stable polygons  (\cite{Hu99a}).

\begin{rem} {\sl How do we choose momentum charges?}
In practice, momentum charges in $\gamma$ may be chosen as follows
by hierarchy. First, list $P$ at the top. Then we choose a general
interior point ${\bf r}$, which serves as the principal charge. On
the next level, we list all coherent subdivisions that
 only refine $P$ but no others. For any polytope $D$ occurring
in the subdivisions of this level, if $D$ contained ${\bf r}$, we
have to stick with it and set $\gamma_S (D) = {\bf r}$. Otherwise,
we may choose freely in the interior of $D$ as long as the
compatibility condition (3) is satisfied. Then we move on to the
third level and list all subdivisions that only refine the
subdivisions from the previous level. Again, for any polytope $D$
in this level, if it contains a charge from the previous level, we
stick with it. Otherwise we choose freely in the interior of $D$
as long as Condition (3) is satisfied. We can go on with this
process until all subdivisions are addressed.
\end{rem}

\subsection{Stable orbits with prescribed  momentum charges}
\label{stableorbits}

\begin{defn}
\label{defn:stableOrbit}
 Fixed a set $\gamma$ of momentum charges.
 A finite  collection of $K$-orbits,
${\bf O} = \{ O_i\}_i$,
is called a stable orbit with momentum charges $\gamma$ if
\begin{enumerate}
\item there is a point $q \in X/\!/^{ch}G$  such that $\cup_i
\overline{G \cdot O_i}$ equals to the support $|F(q)|$; \item for
each polytope $D$ in the subdivision $S= \Phi(|F(q))|$, there is a
unique orbit $O_i$ in $\Phi^{-1} (\gamma_S (D))$. (We will often
denote this orbit by $O_D$.)
\end{enumerate}
In this case, we will say that the stable orbit ${\bf O}$ is of
type $S = \Phi ( |F (q)|)$.
\end{defn}

Observe that each set $G \cdot O_i$ is a single $G$-orbit. Since
$D$ is of top dimensional, and $\gamma_S (D)$ is in the interior
of $D$, it follows from (2) that $\cup_i G \cdot O_i = \cup_i G
\cdot x_i$ where $|F(q)| = \cup_i \overline{G \cdot x_i}$ for some
$x_i \in X$. In fact, if we need, we can even pick $x_i \in O_i$.
For a depiction of stable $K$-orbits, consult Figure 3.

\begin{rem} From the definition, for any stable orbit ${\bf O}$,
there must be an orbit $O_i$ with the principal momentum charge ${\bf r}$.
We will denote it by $O_{\bf r}$, and name it as the principal or main orbit. All  other
orbits will be referred as {\it bubble} orbits of $O_{\bf r}$. Moreover,
for every $D$ and a subdivision  $D= \cup_i D_i$, the orbits $O_{D_i}$ (if any) will be called
 {\it bubbles} of $O_D$ (if any).
\end{rem}

\begin{rem} The definition of
a particular stable orbit ${\bf O}$  with  momentum charges
$\gamma$ utilizes (or depends on) only the values of $\gamma$ on a
single subdivision $\Phi(f^{-1}(q))$. It appears that we may as
well define a stable orbit using only the values of $\gamma$ on
the subpolytopes of $\Phi(f^{-1}(q))$, without referring to the
whole $\gamma$. However, to form a meaningful global moduli space,
various stable orbits must have compatible momentum charges.
Hence we insist to associate ${\bf O}$ with the whole  $\gamma$,
even if it only depends on the values of $\gamma$ on just one
particular subdivision.
\end{rem}

\begin{rem}
In the definition of stable orbits, we may allow some orbits to occur with multiplicities
in the same way as orbit closures may occur in the Chow family. This would be useful if one wishes
to construct and utilizes universal families. Since our approach and application are topological,
the multiplicity issue will be suppressed in this paper.
\end{rem}

\subsection{Local moduli and correspondence varieties}
\label{localmoduli}

Let $P$ be the minimal element in ${\mathcal S}$, and $U_P$ the
set of all stable orbits of type $P$.
This is the local moduli space associated to (the non-subdivision)
$P$.

\begin{defn} In general, given any polytopal subdivision $S \in {\mathcal
S}$, let $Z_S$ be the set of all stable orbits of type $S$.
\end{defn}

This is a subset in
$$\Pi_{D \in S} X_{\gamma_S(D)}.$$
We now describe a neighborhood of $Z_S$,
which is to be an incident variety in a product space.


Let $\widetilde{S} = \bigcup_{S' \le S} S'$. We will  think $\widetilde{S}$
as a collection of subpolytopes. Consider the product space,
$$\Pi_{C \in \widetilde{S}} X_{\gamma_{S'}(C)},$$
where $S'$ is any member of $\widetilde{S}$ that contains $C$.
Note that the expression does not depend on the choice of $S'$,
for, if $S''$ is another one that contains $C$, then by the
compatibility of the set of momentum charges, $\gamma_{S'}(C) =
\gamma_{S''}(C)$. We define an analytic  correspondence variety
$$U_S \subset \Pi_{C \in \widetilde{S}} X_{\gamma_{S'}(C)}$$
as follows.

\begin{defn}
\label{neighborhood}
A point ${\bf O} =\{ O_C\}$ of $ \Pi_{C \in \widetilde{S}} X_{\gamma_{S'}(C)}$
belongs to $U_S$ if both of the following are true:
\begin{enumerate}
\item  there is a unique  $S' \le S$ such that the components
$${\bf O}_{S'} = \{ O_C | C \in S'\}$$
form a stable orbit of type $S'$;
\item  the rest of the components are uniquely determined by
$${\bf O}_{S'} =\{ O_C | C \in S'\}$$ by the means as specified below.
For any $D \in S$ (note that $D$ is of top dimension, cf.
Definition \ref{subdivision}), $D$ is contained in a unique $C \in
S'$ since $S$ refines $S'$.
 In this case we set $$O_D = (G \cdot O_C) \cap \Phi^{-1}(\gamma_{S}(D)).$$
For any other polytope $C'' \in S'' \subset \widetilde{S} \setminus (S \cup S')$,
since $S''$ is refined by $S$,
there must be a polytope $D$ of $S$ such that $D \subset C''$ and
$\gamma_{S''} (C'') \in D$. Hence by Definition \ref{charges} (2),
$\gamma_{S''} (C'') = \gamma_S (D)$. Then, in this case, we simply require
$O_{C''}$ to equal to $O_D$.
\end{enumerate}
Observe that the relation used in this definition is  analytic.
\end{defn}

Recall that  $Z_S$ be the set of all stable orbits of type $S$.
From the above, we see that for any $S' \le S$, there is an
injective map
$$Z_{S'} \hookrightarrow U_S$$
because the components in $Z_{S'}$ completely determine the rest in  $U_S$ as in Definition \ref{neighborhood} (2).
After identifying $Z_{S'}$ with its image in $U_S$, we see that
$$
U_S = \bigcup_{S' \le S} Z_{S'} \subset \Pi_{C \in \widetilde{S}} X_{\gamma_{S'}(C)}$$
From here we immediately have
$$U_{S'} \subset U_S$$
whenever $S' < S$. Now, from Definition \ref{neighborhood}, we see
that the incident relation is analytic and the inclusion $U_{S'}
\hookrightarrow U_S$ is an analytic open embedding. That is, we
have

\begin{prop}
$U_S$ is an analytic subset  of $\Pi_{C \in \widetilde{S}}
X_{\gamma_{S'}(C)}$ (not closed in general). Furthermore, $U_{S'}$
is an open analytic subvariety of $U_S$ whenever $S' < S$.
\end{prop}

\subsection{Global moduli of stable orbits}

Now let ${\bfmit M}_{\gamma}$ be the set of all stale orbits of
type $\gamma$. Then  ${\bfmit M}_{\gamma} = \cup_S U_S$. It
follows from the construction that the complex structures on $U_S$
all agree with each other on the overlaps, and it induces a
Hausdorff topology on ${\bfmit M}_{\gamma}$. That is, ${\bfmit
M}_{\gamma}$ is the inverse limit $\lim_S U_S$ of the system
$\{U_S, U_{S'} \hookrightarrow U_S | S' < S\}$. Furthermore, we
have


\begin{thm}
\label{moduliofstabelKorbits}
The moduli space ${\bfmit M}_{\gamma}$
exists as a separated complex variety, and
is homeomorphic to the Chow quotient $X/\!/^{ch}G$.
\end{thm}

\proof We only need to prove the second statement.

Locally on $U_S$, we define a map
$$\alpha_S: U_S \rightarrow X/\!/^{ch}G$$
as follows. For any point ${\bf O} = \{ O_C\} \in U_S \subset
\Pi_{C \in \widetilde{S}} X_{\gamma_{S'}(C)}$, there is a unique
$S'  \le S$ such that ${\bf O}_{S'} = \{ O_C | C \in S'\}$ is a
stable orbit of type $S'$, and the rest components are uniquely
determined by ${\bf O}_{S'}$. Let $q \in X/\!/^{ch}G$ be the point
in the Chow quotient such that ${\bf O}_{S'} \subset |F(q)|$. By
Definition \ref{defn:stableOrbit} and  Assumption
\ref{uniquesupport}, $q$ is unique. Then we can define
$$\alpha_S ({\bf O}) = q.$$
All those locally defined maps apparently agree with each other on
the overlaps, thus they glue together to give a  globally defined
map
$$\alpha: {\bfmit M}_{\gamma} \rightarrow X/\!/^{ch}G.$$

This map has the inverse
$$\beta:X/\!/^{ch}G \rightarrow {\bfmit M}_{\gamma}$$
by sending a point $q \in X/\!/^{ch}G$ to the stable orbit
$${\bf O} = \{O_D | D \in S = \Phi(f^{-1}(q))\}$$
where $O_D =F (q)\cap \Phi^{-1}(\gamma_S (D))$. (One verifies from the definition that ${\bf O}$ is indeed
a stable orbit of type $S$.)

Note that the Chow quotient, $X/\!/^{ch}G$, can  be stratified according to the subdivision type of $\Phi(f^{-1}(q))$. That is,
$$X/\!/^{ch}G = \bigcup_{S \in {\mathcal S}} Y_S,$$
where $Y_S = \{ q \in X/\!/^{ch}G | \Phi(f^{-1}(q)) =S\}$.
Then $$V_S = \bigcup_{S' \le S} Y_{S'}$$
is an open neighborhood of $Y_S$. The restriction of $\beta$ to $V_S$ has the image in $U_S$.
And, the map
$$\beta |_{V_S} : V_S \rightarrow U_S \subset \Pi_{C \in \widetilde{S}} X_{\gamma_{S'}(C)}$$
is defined component-wise by the projections
$$X/\!/^{ch}G \rightarrow X_{\gamma_{S'}(C)},$$
for all $S \in {\mathcal S}$. Hence $\beta$ is analytic, in
particularly, continuous. Since $\beta$ is a continuous bijection
between two compact Hausdorff spaces, it must be homeomorphism. So
is the inverse map $\alpha$.
\endproof

For a stable orbit ${\bf O} =\cup_i O_i$, we have realized an open
neighborhood $U_{\bf O}$ around it with an incident variety in the
product space $\prod_{D \in S} X_{\gamma_S (D)}$ so that it admits
an induced symplectic form on the smooth locus of the
neighborhood. To put these pieces together to get a global
symplectic or Poisson structure is a task worth pursuing. For, if
we compare with the role of ${\bf r}$ in the symplectic quotient
$\Phi^{-1}({\bf r})/K$, it vividly suggests that, just like what
${\bf r}$ does for $\Phi^{-1}({\bf r})/K$,  $\gamma$ should lead
toward symplectic/Poisson structures on the Chow quotient. This
would make ${\bfmit M}_{\gamma}$ a {\it genuine} symplectic
reduction for the Chow quotient,  adding a new correspondence to
the usual ``GIT=Reduction'' picture. This calls for further
investigation.

\subsection{Blowup along arrangement of subvarieties}

\begin{thm}
Let the notation be as before. Then there is a  holomorphic projection
$${\bfmit M}_\gamma \rightarrow X_{\bf r}, \; {\bf O} \to  O_{\bf r}$$
defined by sending a stable orbit ${\bf O}$ to its principal orbit
$O_{\bf r}$. This map restricts to an isomorphism on the open
subset $X^\circ_{\bf r}$.
\end{thm}
\proof
This follows immediately from the construction of ${\bfmit M}_\gamma$.
\endproof

Obviously, under the enlarged KN correspondence, this map corresponds to the algebraic map
$$X/\!/^{ch}G \rightarrow X_{[{\bf r}]}.$$

Every symplectic reduction $X_{\bf r}$ has a decomposition
$$X_{\bf r} = \bigcup_{D}  M_{{\bf r}, D},$$
where $D$ is a subpolytope of $P$, and an orbit $O \in X_{\bf r}$ belongs to $M_{{\bf r}, D}$ if
$\Phi(\overline{G \cdot O}) = D$. We point out that $M_{{\bf r}, D}= \emptyset$
unless $D$ contains ${\bf r}$. Hence we may write
$$X_{\bf r} = \bigcup_{{\bf r} \in D}  M_{{\bf r}, D}.$$
When $D$ is the whole polytope $P$, $M_{{\bf r}, P}$ is the open subset of
generic points in $X_{\bf r}$.
Set $$N_{{\bf r}, D}= \bigcup_{C \subset D} M_{{\bf r}, C}.$$
This is closed in $X_{\bf r}$.
The complement of $M_{{\bf r}, P}$ is a union  of closed
subvarieties,
$$\bigcup_{D \ne P} N_{{\bf r}, D}.$$
If the group action is quasi-free\footnote{An action is quasi-free if
all the isotropy subgroups are connected. Using the orbifold/stack language,
this assumption may be removed.} and for some general ${\bf r}$,
$$\{N_{{\bf r}, D} | D \ne P\}$$ form an arrangement of smooth subvarieties
(see Definition 1.2 of \cite{Hu99b}), then, we expect that
the projection map $${\bfmit M}_\gamma  \rightarrow X_{\bf r}$$
is a {\it  blowup along the arrangement of smooth subvarieties}, in the sense of
Theorem 1.1 of \cite{Hu99b}. For example, this is the case
for the maximal torus action on the Grassmannian $\Gr(2, {\bf C}^n)$
(Theorem 6.5  \cite{Hu99a} ).

\begin{rem}
The above seems to   provide a (rare) criterion for the smoothness of the Chow quotient,
that is, assuming quasi-free, it is smooth if $\{M_{{\bf r}, D} | D\}$ is an arrangement of smooth subvarieties
of $X_{\bf r}$ for some general ${\bf r}$. This line of approach may be applied to the Chow quotients
of higher Grassmannians $\Gr(n, {\bf C}^m)$ ($n>2, m > 5$), but the combinatorics  involved seems
too intimidating at the moment.
\end{rem}

\bigskip
\section{Chow Quotient: topological approach}

There is yet another topological approach to the Chow quotient,
which is somewhat ``orthogonal'' to the approach of stable $K$-orbits.

\subsection{Action-manifolds}
\label{AM}

Instead of $K$-orbits,
 we can also consider the following infinitesimal action. For any $\xi \in {\bfmit k}$, set
$$\sqrt{-1} \xi_{X, x} := {d \over dt}_{|_{t=0}}
exp ( t \sqrt{-1} \xi) \cdot x.$$
Treating  $\sqrt{-1} {\bfmit k}$ as a distribution of vector fields on $X$,
we obtain its integral manifolds through points of $X$.  In this case, the integral
manifolds are not closed. Hence
we take the closures of these integral manifolds. We will see that the closures are homeomorphic,
via the moment map, to subpolytopes of $\Phi(X)$,  and hence are manifolds with corners, in general.
Thus we may call the above integral manifolds  {\it open  action-manifolds}
and their closures {\it action-manifolds} (with corners),
because they come from the group action.

Let $G = K \cdot A$ be the polar decomposition. Then it can be verified that

\begin{prop} (\cite{Atiyah}, \cite{GS}.)
The open  action-manifold through a point $x$ is $A \cdot x$.
The action-manifold through the point $x$ is $\overline{A \cdot x}$. Moreover,
the moment map $\Phi$ induces a homeomorphism
between $\overline{A \cdot x}$ and $\Phi(\overline{A \cdot x})$.
\end{prop}

For action-manifolds through {\it generic points}\footnote{e.g.,
points of the  open subset $U$ in Definition
\ref{defnChowQuotient}.}, we will call them generic
action-manifolds. Two generic action-manifolds are equivalent if
one can be obtained from the other by the action of an element of
$K$. Let ${\bfmit M}^\circ$ be the moduli space of equivalence
classes of generic action-manifolds. Then ${\bfmit M}^\circ$
contains $U/G$ as an open subset and is itself   an open variety.
We would like to describe a geometrically meaningful
compactification ${\bfmit M}$ of ${\bfmit M}^\circ$ by  providing
natural geometric meanings of the boundary points. These boundary
points will be called stable action-manifolds. Generic
action-manifolds, as generic points of ${\bfmit M}$, are
automatically stable. So, what are the rest stable
action-manifolds? To answer this question, we need some
preparation.

\subsection{Perturb, translate and specialize: topological version}
\label{section:toppush}

Let $x$ and $y$ be two points with $$\dim A \cdot x = \dim A \cdot y =\dim A.$$
Let ${\bf r} = \Phi(x)$.
Take a real analytic slice
\footnote{topological slice will suffice. cf. Remark \ref{top=ana} and also Remark \ref{topnature}.},
$$R_x \subset \Phi^{-1}({\bf r}),$$ around the point $x$, transversal to $K$-orbits.

Recall that a point of $X$ is said to be isotropy-free if its isotropy subgroup is the identity.

\begin{defn}
\label{toppush}
We say $x$ can be  perturbed (to general positions), translated (along $A$-orbits), and specialized
to $y$, which is
denoted by $$x  \to_A y,$$
if there is a generic real holomorphic map from the interval $I=[-1,1]$ to the slice $R_x$,
$$\varphi: I \rightarrow R_x,$$
with $\varphi (0) = x$ such that
 $\varphi (t)$ is a generic isotropy-free point for all $t \ne 0$,
and in addition there is a  real holomorphic map from the punctured interval
$I^*=I \setminus\{0\}$ to the group $A$
 $$a: I^* \rightarrow A$$
such that $\psi(t)= a(t) \cdot \varphi (t) \in \Phi^{-1}(\Phi(y))$,  and
$$y=\lim_{t \to 0} a(t) \cdot \varphi (t).$$
\end{defn}

Just as the original P.T.S, the above relation is symmetric. To prove this, we can choose a real analytic slice
$R_y$ containing $\psi(t)= a(t) \cdot \varphi (t)$, and then repeat the arguments of Remark \ref{symetric}.
 That is, if $x  \to_A y$, then $y  \to_A x$. Because of this, we may write $ x \sim_A y$.

\begin{rem}
\label{top=ana} In the above definition, we can replace ``real
holomorphic map'' by ``continuous map'', all the statements and
proofs, which are all topological in nature,
 remain unchanged.
\end{rem}

In Definition \ref{toppush}, after the choice of the map $\varphi$ is made, $a(t)$ is  uniquely determined.

\begin{lem}
\label{lem:a(t)unique}
For  $t \ne 0$,
$a(t)$ is the unique point $a$ in $A$ such that
$\Phi (a \cdot \varphi(t)) = \Phi(y)$.
\end{lem}
\proof
For any  $t \ne 0$, since $\varphi(t)$ is isotropy-free and generic,
 we have a homeomorphism
 \begin{equation*}
\begin{CD}
A @>{\cong}>>  A \cdot \varphi(t) @>{\Phi}>> P^\circ \end{CD}
\end{equation*}
 where $P^\circ$ is the interior of $P$. Hence $a(t)$ is the
unique point $a$ in $A$ such that $\Phi (a \cdot \varphi(t)) =
\Phi (y) \in P^\circ$.
\endproof

\subsection{Stable action-manifolds}

\begin{defn}
\label{defn:sAM}
A finite  union of action-manifolds, $\cup_i \overline{A \cdot x_i}$,
is called a stable action-manifold if
\begin{enumerate}
\item $\dim A \cdot x_i = \dim A$ for all $i$;
\item $x_i \sim_A x_j$ for all $i$ and $j$;
\item the moment map $\Phi$ induces a homeomorphism
between $\cup_i \overline{A \cdot x_i}$ and $\Phi(X)$.
\end{enumerate}
\end{defn}

Note that the condition (3) implies that $\Phi(\cup_i \overline{A
\cdot x_i}) = \Phi(X)$. It also implies that $\cup_i \overline{A
\cdot x_i}$ is connected subset of $X$. In particular, the indexes
can be re-arranged so that $\overline{A \cdot x_i} \cap
\overline{A \cdot x_{i+1}} \ne \emptyset$. (Consult Figure 4 for
an illustration.)

\begin{prop}
\label{inchowfiber}
If $x \sim_A y$, then $x$ and $y$ are in the same Chow fiber. In particular,
every given stable action manifold is contained in a single Chow fiber.
\end{prop}
\proof
After modulo the action of $K$, we can treat orbit $A \cdot x$ as point $[G \cdot x]$
in GIT quotient. Note also that the proof of Theorem \ref{chowfiberAG} only use the continuity of the maps
involved, but not the holomorphic properties (cf. Remark \ref{topnature}).
Henceforth the proofs of  Theorem \ref{chowfiberAG} (the sufficient part) can  be repeated almost
word by word to conclude the statement of this proposition. Further details are omitted.
\endproof

Thus if a stable action-manifold, ${\bf M}= \cup_i \overline{A \cdot x_i}$, is contained in a Chow fiber $f^{-1}(q)$, then
$S = \cup_i \Phi(\overline{A \cdot x_i}) = \Phi (f^{-1}(q))$ is a subdivision of $P= \Phi(X)$. In this case,
we will say that ${\bf M}$ corresponds to the subdivision $S$ or is of type $S$.

\begin{lem}
\label{k1=k2}
Suppose that we have $x \sim_A x'$ and $y \sim_A y'$. Assume further that
$y =k_1 \cdot x$ and
$y' = k_2 \cdot x'$ for some $k_1, k_2 \in K$.
Then $k_1=k_2$.
\end{lem}
\proof
Assume that we have $$\lim_{t \to 0} a(t) \cdot \varphi (t) = x'$$ with $\varphi (0) = x$.
Then $\psi (t) = k_1 \cdot \varphi (t)$ defines a generic real holomorphic map from
the interval $I$
to a slice $R_{y}$ near $y$ with $\psi (0) =k_1 \cdot x = y$.
Let $b(t)$ be as in Definition \ref{toppush} such that
$$\lim_{t \to 0} b(t) \cdot \psi (t) = y'.$$
Then  because $\psi (t) = k_1 \cdot \varphi (t)$ and $y' = k_2 \cdot x'$,  we obtain
$$\lim_{t \to 0} b(t) k_2^{-1}k_1 \cdot \varphi (t) = x'.$$
By Lemma \ref{lem:a(t)unique}, for generic $t$, we must have
$b(t) k_2^{-1}k_1 = a(t)$, that is,
$$b(t)a(t)^{-1} = k_2 k_1^{-1}.$$
Since $A \cap K = \id$, we have $k_1 = k_2$.
\endproof

\begin{defn}
Two stable action-manifolds are said to be equivalent if there is an element $k \in K$ such that
the action of $k$ sends one stable action-manifold to the other.
\end{defn}

It is immediate that two equivalent action manifolds must lie in
the same Chow fiber because any Chow fiber is $G$-invariant. The
following proves the converse.

\begin{prop}
\label{sameChowFiber} If two stable action manifolds ${\bf M}_1$
and ${\bf M}_2$  are in the same Chow fiber, then they are
equivalent.
\end{prop}
\proof Assume that ${\bf M}_1= \cup_i \overline{A \cdot x_i}$ and
${\bf M}_2= \cup_i \overline{A \cdot y_i}$ are in the same Chow
fiber. Then we have $\cup_i \overline{G \cdot x_i} = \cup_i
\overline{G \cdot y_i}$. By Definition \ref{defn:sAM} and
Proposition \ref{inchowfiber}, we can re-arrange so that there is
a one-to-one correspondence between $\{x_i\}$ and $\{y_i\}$, and
$G \cdot x_i = G \cdot y_i$. By choosing different representatives
of $A \cdot x_i$ (for all $i$) if necessary, we may assume that
$$y_i = k_i \cdot x_i$$
for some $k_i \in K$ for all $i$.
Now apply Lemma \ref{k1=k2}.
\endproof

\begin{rem}
This proposition shows that
a stable action-manifold is not just an arbitrary union $\cup_i \overline{A \cdot x_i}$
of $A$-orbits even if we require that $\cup_i \overline{G \cdot x_i}$ is a Chow fiber.
Geometrically, stable action-manifolds occur as the
limiting configurations of families of generic $\overline{A \cdot x}$.
\end{rem}



\subsection{Moduli of stable action-manifolds }
We will use ${\bfmit M}$ to denote the set of equivalence classes of all stable action-manifolds.
(Note that the definition of ${\bfmit M}$ involves no choices; while the moduli space of stable
$K$-orbits with momentum charges $\gamma$, as it depends on $\gamma$, is always denoted by ${\bfmit M}_\gamma$.)

\begin{thm}
\label{MofSA} The moduli space ${\bfmit M}^\circ$ of generic
stable action-manifolds admits a natural compactification ${\bfmit
M}$ by adding stable action-manifolds. The resulting space is
analytic, and is homeomorphic to the Chow quotient $X/\!/^{ch}G$.
\end{thm}

\proof
The approach to this theorem, although somewhat ``orthogonal''
to that of Theorem \ref{moduliofstabelKorbits},
is in spirit related and similar to it. We will only give a sketch.

Let $G = A \cdot K$ be the polar decomposition. First, recall that
every piece $M$ in a stable action-manifold ${\bf M}$ is of the
form $\overline{A \cdot x}$. Since $K$ is compact, there is
one-to-one correspondence between (the closures of) $G$-orbits
and (the closures of) $A$-orbits modulo $K$, and hence we may
identify the two kinds of orbits in GIT quotient  $X^{ss}/\!/G$.
In other words, we may write $[A \cdot x]$ for $[G \cdot x]$ in
$X^{ss} /\!/G$.

The moduli space ${\bfmit M}$ is canonically defined (depends on no choices).  However, to prove this theorem,
 we have to make some auxiliary choices. That is,
we will  fix a set $\gamma$ of momentum charges. By a small perturbation, we may require
that all the charges are rational. (This technical maneuver is only needed to allow us
to use GIT quotients.)

By Proposition \ref{inchowfiber},
any  stable action manifold ${\bf M}'$ is contained in some Chow fiber $f^{-1}(q)$.
Hence it corresponds to some subdivision $S = \Phi( f^{-1}(q))$. Using the convention mentioned in the
beginning of the proof, we will embed ${\bf M}'$ in the product space of some GIT quotients,
$$\Pi_{C \in  \widetilde{S}} X_{[\gamma_ {S'} (C)]},$$
and then define an open neighborhood $W_S$ of ${\bf M}'$ as an
incident analytic subvariety in $\Pi_{C \in  \widetilde{S}}
X_{[\gamma_ {S'} (C)]}$. Here  $\widetilde{S} = \bigcup_{S' \le S}
S'$, and
 $S'$ is any member of $\widetilde{S}$ that contains $C$. The product space does not depend on the choice of $S'$
(cf. the remark in the paragraph immediately before Definition
\ref{neighborhood}). As remarked in the beginning of the proof, we
will represent a point of $X_{[\gamma_ {S'} (C)]}$ by an $A$-orbit
closure.

A point ${\bf M} =\{ [M_C]\}$ of $\Pi_{C \in \widetilde{S}} X_{[\gamma_{S'}(C)]}$
belongs to $W_S$ if both of the following are true:
\begin{enumerate}
\item  there is a unique  $S' \le S$ such that the components
$${\bf M}_{S'} = \{ [M_C] | C \in S'\}$$
is a stable action-manifold corresponding to the subdivision $S'$;
\item  the rest of the components are completely determined by ${\bf M}_{S'}$ as follows.
For any $D \in S$,  $D$ is contained in a unique $C \in S'$ since $S$ refines $S'$.
In particular, $C$ contains the Chamber $[\gamma_S(D)] \subset D$.
In this case, we require $$[M_D] = [M_C] \in X_{[\gamma_S(D)]}.$$
Here, using the remark in the beginning of the proof, we may treat $M_C$, originally
an orbit (closure) of type $C$ (i.e., $\Phi (M_C) =C$), as an orbit in $X_{[\gamma_S(D)]}$ as well.
For any other polytope $C'' \in S'' \subset \widetilde{S} \setminus (S \cup S')$,
since $S''$ is refined by $S$,
there must be a polytope $D$ of $S$ such that $\gamma_{S''} (C'') = \gamma_S (D)$
(cf. Definition \ref{neighborhood} (2)), and in this case, we require that
$[M_{C''}]$  equals $[M_D]$.
\end{enumerate}
This makes $W_S$ an analytic subvariety of $\Pi_{C \in \widetilde{S}} X_{[\gamma_{S'}(C)]}$.
As in the case of  stable $K$-orbits, after the obvious identifications,
we have that $W_{S'} \subset W_S$ whenever $S' < S$.
In particular,  all these complex structures
agree with each other on the overlaps.
Consequently, we obtain that
the moduli space ${\bfmit M}$
is a separated complex analytic variety.

Now using basically the same argument as in the proof of Theorem \ref{moduliofstabelKorbits},
we can define a map
$$\theta: {\bfmit M} \rightarrow X/\!/^{ch}G$$ and its inverse
$$\theta^{-1}: X/\!/^{ch}G \rightarrow {\bfmit M},$$
and prove that ${\bfmit M}$ is homeomorphic to the Chow quotient $X/\!/^{ch}G$.
Further details are omitted.
\endproof




\bigskip\bigskip



\end{document}